\documentclass[11pt,reqno]{article}

\usepackage{latexsym}
\usepackage{amssymb}

\def\ZZ {{\mathbb Z}}
\def\NN {{\mathbb N}}

\def\RR {{\mathbb R}}

\def\Si{\Sigma}
\def\La{\Lambda}
\def\De{\Delta}
\def\Om{\Omega}
\def\Ga{\Gamma}

\def\Ups{\Upsilon}

\def\la{\lambda}
\def\ka{\kappa}
\def\be{\beta}
\def\de{\delta}

\def\ga{\gamma}
\def\ve{\varepsilon}

\def\pfrk{\pitchfork}

  \def\cG{{\cal G}}  
  \def\cH{{\cal H}}  \def\cT{{\cal T}}
\def\cC{{\cal C}}  \def\cI{{\cal I}}

   \def\cR{{\cal R}}

\def\mfT{{\mathfrak{T}}}
\def\mfG{{\mathfrak{G}}}
\def\mfR{{\mathfrak{R}}}

\newtheorem{theor}{Theorem}[section]
\newtheorem{lemm}[theor]{Lemma}
\newtheorem{clai}[theor]{Claim}
\newtheorem{bdl}[theor]{Bounded Distortion Lemma}
\newtheorem{fact}[theor]{Fact}
\newtheorem{coro}[theor]{Corollary}

\newtheorem{algo}[theor]{Algorithm}
\newtheorem{prop}[theor]{Proposition}

\newtheorem{rema}[theor]{Remark}

\newtheorem{maintheorem}{Theorem}

\newenvironment{demo}{
\noindent{\bf Proof: }}{\hfill$\Box$\medskip}

\title{Collision, explosion and collapse of homoclinic classes}
\author{Lorenzo J. D\'\i az and Bianca Santoro
\thanks{This paper was partially supported by CAPES, CNPq, Faperj,
and Pronex Dynamical Systems (Brazil)}}

\begin{document}

\date{July 4th,  2004}
\maketitle

\begin{abstract}
Homoclinic classes of generic $C^1$-diffeomorphisms are
maximal transitive sets and pairwise disjoint.
We here present a model explaining how two different homoclinic classes
may intersect, failing to be disjoint.
For that we construct a one-parameter family of diffeomorphisms
$(g_s)_{s\in [-1,1]}$ with hyperbolic points $P$ and $Q$
having nontrivial homoclinic classes,
such that, for $s>0$,
the classes of $P$ and $Q$ are disjoint,
for $s<0$, they are equal, and, for $s=0$, their intersection
is a saddle-node.
\end{abstract}

\section*{Introduction}

In this paper we study the collision of 
non-trivial homoclinic classes
{\em via\/} saddle-node bifurcations and the dynamics before and
after this collision. 
The main motivation of this paper comes from recent results 
about {\em maximal transitive sets\/}:
for 
generic\footnote{by a generic diffeomorphism we mean a 
diffeomorphism in a 
residual subset $\cR$ of Diff$^1(M)$.}
$C^1$-diffeomorphisms, the homoclinic classes 
are either disjoint or coincide, see \cite{Ar} and \cite{CMP}.

Let us start by recalling some basic definitions.
Given a diffeomorphism $f$,
 an $f$-invariant set
$\La$ is {\em transitive\/} 
if there is an $x \in \La$ whose forward orbit is dense in $\La$, i.e.,
$\La = \overline{\cup_{i \in \NN}f^i(x)}$.
A transitive set is {\em maximal\/} if it is a maximal element
of the family of all transitive sets partially ordered by inclusion.
Observe that any transitive set is contained in a maximal one.
A transitive set $\La$ is {\em saturated\/} if it 
contains every transitive set $\Si$
such that $\La\cap \Si \ne\emptyset$.
Clearly, every saturated transitive set is also maximal.
The {\em homoclinic class\/} of a saddle $P$ of
$f$, denoted by $H(P,f)$, is the closure of the transverse intersections 
of the orbits of the stable and unstable manifolds of
$P$. Every homoclinic class is a transitive set, not necessarily
maximal nor saturated.

The problem of characterizing  and
describing (for a large class of systems)
maximal and saturated transitive sets is a key problem in dynamics.
In fact, these saturated transitive sets are the natural candidates for
playing the role of the elementary pieces of dynamics
(similar to the role of the basic sets in the 
Smale hyperbolic theory, \cite{Sm}).
Recently,  \cite{Ab}  states that for generic 
$C^1$-diffeomorphisms $f$ having finitely many different homoclinic classes
the non-wandering set of $f$, $\Om(f)$, is the disjoint union
of such classes. Moreover, these  classes 
verify a weak form of hyperbolicity (existence of a dominated 
splitting, see \cite{BDP})
and are the maximal invariant sets
of a fixed filtration (see Section~\ref{ss.endofthmain})
independent of the 
generic diffeomorphism in 
a neighborhood of $f$.

Consider a closed manifold $M$ and
denote by Diff$^1(M)$ the space of $C^1$-diffeomorphisms
endowed with the usual uniform topology. 
In \cite{Ar}, it is proved that homoclinic classes of  
generic
diffeomorphisms are {\em maximal transitive sets.\/}
\cite{CMP} generalizes this result by proving that homoclinic
classes of generic diffeomorphisms are saturated
transitive sets.
Thus  homoclinic classes
of generic diffeomorphisms are either equal or disjoint.
Let us observe that there are locally generic diffeomorphisms
having saturated transitive sets without periodic orbits (so
which are not homoclinic classes), see \cite{BD01}.

The goal of this paper is to give examples of
homoclinic classes which are 
not saturated transitive sets,  
presenting an
explanation of how this pathology arises.
In fact, we exhibit homoclinic classes
which are not contained in any saturated transitive set.
For simplicity, we  consider diffeomorphisms
defined on three manifolds, but our constructions can be
carried out to higher dimensions after
straightforward modifications.

The homoclinic  classes 
in this paper simultaneously contain
(in a stable way) 
hyperbolic fixed  points 
of {\em Morse index\/} (dimension of the unstable bundle)
one and two, hence not hyperbolic.
We construct a diffeomorphism $f$ with 
saddles $P$ and $Q$ of Morse indices one and two,
such that 
their homoclinic classes  are  nontrivial,
maximal transitive,
whose intersection is just a saddle-node (so the classes are not saturated
transitive sets).

\begin{maintheorem}
\label{t.main}
There exist an open set $W$ and  a family of diffeomorphisms
$(g_s)_{s\in [-1,1]}$ such that, for every $s$, the diffeomorphism
$g_s$ has hyperbolic fixed points $P$ and $Q$ of Morse indices
$1$ and $2$ such that the maximal invariant set of $g_s$ in $W$,
denoted by $\La_s$, verifies the following:
\begin{itemize}
\item For every small $s<0$, 
the set $\La_s\cap \Om(g_s)$ 
is the
disjoint union of $H(P,g_s)$ and $H(Q,g_s)$,
where $H(P,g_s)$ and $H(Q,g_s)$ are non-hyperbolic and locally maximal. 
\item 
 For $s=0$, $\La_s  = H(P,g_s) \cup H(Q,g_s)$,
where $H(P,g_s)$ and $H(Q,g_s)$ are locally maximal
and  
 $H(P,g_s)\cap H(Q,g_s)=\{S\}$, where $S$ is a saddle-node fixed point.
\item
 For every small $s>0$, 
$\La_s=H(P,g_s)=H(Q,g_s)$.
\end{itemize}
\end{maintheorem}

This result means that the homoclinic classes of $P$ and $Q$
{\em collide\/} at $s=0$ and thereafter {\em explode\/} (the point $P$ that
does not belong to $H(Q,g_0)$ is in 
$H(Q,g_s)$ for positive $s$, and the same holds for
the point $Q$ and $H(P,g_s)$). Finally, these two homoclinic classes
also {\em collapse\/} 
($H(Q,g_s)=H(P,g_s)$) for positive $s$.

Taking the set  $W$ to be a level of a 
{\em filtration,\/} (see Section~\ref{ss.endofthmain}), 
one gets the following:

\begin{maintheorem}
\label{t.saturated}
The homoclinic classes 
$H(P,g_0)$ and $H(Q,g_0)$ are not saturated and 
they are not contained
in any saturated transitive set. 
\end{maintheorem}

\medskip

Our construction involves 
{\em saddle-node bifurcations\/} and
{\em heterodimensional cycles.\/}
In fact, we introduce a codimension-two bifurcation, the 
{\em saddle-node heterodimensional cycles,\/}
and study the {\em lateral homoclinic classes of a saddle-node.\/}
Let us explain all that in details.

Consider a diffeomorphism $f$ having two hyperbolic fixed points
$P$ and  $Q$ with Morse indices $1$ and $2$, 
respectively. Then, 
$f$ has a  
{\em heterodimensional cycle  associated to $P$ and $Q$\/}
if the stable manifold of  $P$, denoted by $W^s(P,f)$,
and the unstable manifold of $Q$, $W^u(Q,f)$,
have a non-empty transverse intersection,
and 
the unstable manifold of $P$, $W^u(P,f)$,
and the stable one of  $Q$, $W^s(Q,f)$, 
have a quasi-transverse intersection throughout the orbit of a point
$x_0$, i.e., $T_{x_0} W^s(Q,f) + T_{x_0} W^u(P,f) 
= T_{x_0} W^s(Q,f) \oplus T_{x_0} W^u(P,f)$,  thus
 $\dim(T_{x_0} W^s(Q,f) + T_{x_0} W^u(P,f))=2$.
Bifurcations through heterodimensional cycles have been
sistematically studied in the series of papers
\cite{Detds,Dnon, DRnc, DRnon, DU, DRetds,DRfund}.  

A {\em saddle-node} $S$ of a diffeomorphism
$f$ is a periodic point (we here assume to be fixed)
such that the derivative of $f$ at $S$ has 
$1$ as its only eigenvalue in the unitary circle.
We consider
saddle-nodes of saddle type (i.e., the derivative of $f$
at $S$ simultaneously  has eigenvalues inside and outside
the
unitary circle). 
This means that the tangent bundle of $M$ at $S$ 
has a $Df$-invariant splitting
$E^{ss}\oplus E^c \oplus E^{uu}$,
where $E^{ss}$ (resp. $E^{uu}$) is the bundle spanned by the eigenvectors 
associated to the contracting (resp. expanding) eigenvalues, and $E^c$
is the eigenspace associated to the eigenvalue $1$
(in our context, all these spaces have dimension $1$).
According to the theory of invariant manifolds,
see \cite{HPS}, there are defined the {\em strong stable\/}
and {\em unstable\/} manifolds of the saddle-node,
defined as the unique $f$ invariant manifolds tangent at $S$  to
$E^{ss}$ and to $E^{uu}$ and denoted by $W^{ss}(S,f)$
and $W^{uu}(S,f)$, respectively.

Motivated by the fact  that (generic) saddle-nodes 
(of saddle type) simultaneously behave  as 
 points of index two and 
one (the stable and unstable manifolds of the saddle-node
have both dimension $2$), 
we define
{\em saddle-node heterodimensional cycles.\/}
A diffeomorphism $f$ has a saddle-node heterodimensional cycle 
associated to a saddle-node $S$ and the saddle $P$
of Morse index one if the (two-dimensional) unstable manifold of $S$
and stable manifold of $P$ have nonempty
transverse intersection and the (one-dimensional) invariant manifolds
 $W^{ss}(S,f)$ and $W^{u}(P,f)$
have a quasi-transverse intersection along the orbit of
some point.
One similarly  defines saddle-node heterodimensional cycles 
associated to a saddle-node $S$ and a saddle $Q$
of Morse index two. 

Roughly speaking, in our construction we  consider a diffeomorphism
$f$  simultaneously having two saddle-node heterodimensional cycles.
We consider a two parameter family $(f_{t,s})_{t,s\in [-1,1]}$ 
of diffeomorphisms
such that $f_{0,0}$ has a pair of saddle-node heterodimensional cycles,
one associated to
a saddle-node $S$ and a saddle $P$  of Morse index one and other one
associated to a saddle $Q$ of index one
and the saddle-node $S$. The parameter $t$ describes the unfolding of the  
cycles (relative motion between compact parts of $W^u(P,f_{t,0})$ and
$W^{ss}(S,f_{t,0})$
and of $W^s(Q,f_{t,0})$ and
$W^{uu}(S,f_{t,0})$).
The parameter $s$ describes the unfolding of the saddle-node:
for positive $s$ there are two saddles $S^+_s$ and $S_s^-$ of
indices $2$ and $1$, colliding at $s=0$ to the saddle-node $S$
and disappearing for negative $s$.
We see that, fixed a small $\bar t>0$,
for $s>0$
(before the collapse of the saddles),
$H(P,f_{\bar t,s})=
H(S_s^+,f_{\bar t,s})$ and 
$H(Q,f_{\bar t,s})=
H(S_s^-,f_{\bar t,s})$ for all small positive $s$. Moreover,
 $H(P,f_{\bar t,s})\cap H(Q,f_{\bar t,s})=\emptyset$.
At the saddle-node bifurcation we have 
$H(P,f_{\bar t,0})\cap H(Q,f_{\bar t,0})=\{S\}$. Finally,
for $s<0$,
after the disappearing of the saddles, 
$H(P,f_{\bar t,s})=H(Q,f_{\bar t,s})$.
See the results in Section~\ref{s.colapso}.
Theorem~\ref{t.main} follows by considering the arc
$g_s=f_{\bar t,-s}$.
To deduce Theorem~\ref{t.saturated} from Theorem~\ref{t.main},
we consider a filtration
having the open set $W$ as a level and 
analyze the orbits of recurrent points of $\La_s$.

In forthcoming papers, we will illustrate
how this type of bifurcation naturally appear 
as secondary bifurcations in the unfolding of heterodimensional
cycles and give the model for the collision, 
explosion, and collapse of (nontrivial) hyperbolic homoclinic classes,
see \cite{DRprep}.

\medskip

Let us say a few words about our constructions.
As mentioned, our setting necessarily corresponds
to a non-generic situation, so we focus our attention on an
example (we have not done any effort for generality).
We begin by presenting (in Section~\ref{s.heterodimensional})
 a model for the unfolding 
of a heterodimensional cycle. This model (motivated
by \cite{Detds}, \cite{BD1}, \cite{DRetds} and \cite{DRfund}) 
allows us to give
a rather transparent explanation of the dynamics in the unfolding
of a cycle by
reducing it to the study of the dynamics of an iterated system
of functions defined on an interval, this is done in
Section~\ref{s.dinamica1dim}. 
Recall that the dynamics of a (linear) Smale horseshoe
is given by two affine expanding maps of the interval
(say $I=[0,1]$)
whose domains of definition are two disjoint closed subintervals of
$I$ (say $[0,1/3]$ and $[1/3,1]$).
The interval $(1/3,2/3)$ 
is the main gap of the horseshoe and corresponds to points in the basin
of attraction of a sink.
The affine model 
associated to  heterodimensional cycles  
is given by infinitely many expanding affine maps 
$F_i$ defined on subintervals $I_i$ of $I$ 
which are non-disjoint
(the interior of the intervals
$I_i$ are pairwise disjoint, but $I_i$ and $I_{i+1}$
have a common extreme). 
Thus in this model there are no gaps and
there are no escaping points.

In Section~\ref{s.pontosnaoerrantes}, we prove that,
after  unfolding of
the cycle, 
the dynamics of  the model family is non-hyperbolic: the 
point of index $1$ in the cycle belongs
to  the homoclinic class of the point of index $2$
in the cycle. In fact, in this section, using the one-dimensional reduction,
we give a
shorter and clearer proof of the results in
\cite{Detds}.  
Since all the constructions in this paper
rely heavily on this proof and there is not any
written version of this approach, we have decided to include 
a short description of them.

In Section~\ref{s.selano},
we introduce the 
{\em lateral homoclinic classes of a saddle-node $S$\/}
of saddle type of a diffeomorphism $f$ as above,
$H^+(S,f)$ and $H^-(S,f)$, respectively 
defined as the closure of the transverse intersections 
$W^u(S,f)\pfrk W^{ss}(S,f)$ and
$W^s(S,f)\pfrk W^{uu}(S,f)$. 
These lateral homoclinic classes essentially behave as the
usual ones.
We see that for arcs $f_t$ unfolding at $t=0$ the  saddle-node
heterodimensional cycle (associated to the saddle of index one $P$ and
the saddle-node $S$)
one has 
$H(P,f_t)\subset H^+(S,f_t)$ for all small positive $t$.
Moreover, under mild conditions, one gets 
$H^+(S,f_t)=H(P,f_t)$ for all small $t>0$. 
The inclusion $H(P,f_t)\subset H^+(S,f_t)$ follows adapting 
(in a rather
straightforward way)
the results for the model family in Section~\ref{s.pontosnaoerrantes}. 
For the inclusion
$H^+(S,f_t)\subset H(P,f_t)$ we need new ingredients that we
borrow from \cite{Dnon}. 

Using the
results in Sections~\ref{s.pontosnaoerrantes}
and \ref{s.selano}, we get a complete description of
the homoclinic classes $H(P,f_{\bar t,s})$ and  
$H(Q,f_{\bar t,s})$ before the collapse of the saddles
$S^+_s$ and $S^-_s$ to the saddle-node.
Finally, to study
$H(P,f_{\bar t,s})$ and  
$H(Q,f_{\bar t,s})$
after the collision,
we introduce new systems of iterated
(one dimensional) functions and analyze their dynamics.

\medskip

\noindent
{\bf Acknowledgements:}
We thank W. Horita for the careful reading of a first
version of this paper and  C. Bonatti for suggesting 
the use of the model family in Section~\ref{s.heterodimensional}.

\medskip

\section{Heterodimensional cycles: a model family}
\label{s.heterodimensional}

In this section, we construct a model one-parameter  
family $(f_t)_{t\in [-1,1]}$ 
of diffeomorphisms unfolding a heterodimensional cycle.
The study of the semi-local dynamics of $f_t$
will be reduced to the analysis of a one-parameter family of 
endomorphisms with infinitely many discontinuities which
describe the dynamics of $f_t$ in the {\em central direction,\/}
see Section~\ref{s.dinamica1dim}.

Consider a diffeomorphism $f$ with a 
heterodimensional cycle having the  following dynamical
configuration. In local coordinates in $\RR^3$, the cycle is
associated to saddle fixed  points
$Q=(0,0,0)$ and
$P=(0,1,0)$ of indices $2$ and $1$, respectively, 
verifying the
following conditions:

\medskip
\noindent
{\bf Partially hyperbolic (semi-local) dynamics of the cycle:}
\begin{itemize}
\item
In the set
$[-1,1]\times [-1,2]\times [-1,1]$ the diffeomorphism has the 
form
$$
f(x,y,z)=
(\la_s x, F(y),
\la_u z),
$$
where
$F\colon [-1,2]\to (-1,2)$ is a strictly increasing
 monotone map with exactly two
fixed points, a source at $0$ and a sink at $1$,  
and
$0<\la_s < d_m < d_M  < \la_u$,  where
$0< d_m
<F^\prime(x)<d_M$
for all $x\in [0,1]$.
\item
 There is $\delta>0$ such that
$F$ in linear in $[-\delta,\delta]$ and 
affine in $[1-\delta,1+\delta]$.
We denote by $\beta>1$ and $0<\lambda<1$,
the eigenvalues of $F$ at
$0$ and $1$, respectively.
\end{itemize}
Observe that
$[-1,1]\times \{(0,0)\}\subset W^s(Q)$,
$\{0\}\times [0,1)\times [-1,1]
\subset W^u(Q)$,
$\{(0,0)\}\times   [-1,1]\subset W^u(P)$,
and
$[-1 ,1]\times (0,1]
\times \{0\}
\subset W^s(Q)$.
Thus 
$\ga= \{0\}\times (0,1)\times \{0\}$ is a normally hyperbolic
curve contained in $W^u(Q)\cap W^s(P)$.

\medskip
\noindent
{\bf Existence and unfolding of the cycle:}
\begin{itemize}
\item
{\bf The cycle:}
There exist $k_0\in \NN$ 
 and a small neighborhood  $U$ of 
$(0,1,-1/2)\in W^u(P)$ 
such that the restriction of $f^{k_0}$
to $U$ is a translation,
$$
f^{k_0}(x,y,z)=
(x-1/2, y-1, z+1/2).
$$
In particular,
$f^{k_0}(0,1,-1/2)=
(-1/2,0,0)\in W^u(Q)$. Thus 
$W^s(Q)$ and $W^u(P)$ meet
throughout the orbit of the heteroclinic point $(-1/2,0,0)$.
By construction,
$(-1/2,0,0)$ is a quasi-transverse heteroclinic point.
\item
{\bf The unfolding of the cycle:} 
Consider the one parameter family $(f_t)$
of diffeomorphisms  coinciding 
with $f$ in 
$[-1,1]\times [-1,2]\times [-1,1]$ and 
such that the restriction of $f_t^{k_0}$ to $U$ is of the
form
$$
f^{k_0}_t(x,y,z)=
(x-1/2, y-1+t, z+1/2)
=
f^{k_0}(x,y,z)+(0,t,0).
$$
Therefore, for small $t>0$,
$\{(-1/2,t)\}\times [-1,1]\subset W^u(P,f_t)$. Thus
$x_t=(-1/2,t,0)$ is a transverse homoclinic point of
$P$ (for $f_t$). Similarly, $y_t=(-1/2, 0,0)$ is a transverse homoclinic
point of $Q$ (for $f_t$).
\end{itemize}

Consider a small {\em neighborhood of the heterodimensional cycle
associated with $f_0$,\/} that is, an open set $W$ containing the
connexion curve $\ga=\{0\}\times [0,1]\times \{0\}$
and the $f_0$-orbit of the heteroclinic point
$(-1/2,0,0)$.
For small $t$, let $\La_t$ 
be the {\em maximal $f_t$-invariant set in $W$,\/}
$\La_t = \cap_{n\in\ZZ}f_t^{n}(W)$.
We also consider the forward and backward invariant sets in $W$, 
$\La_t^+ = \cap_{n\ge 0}f_t^{-n}(W)$ and
$\La_t^- = \cap_{n\ge 0}f_t^{n}(W)$.

Fix a small positive $\rho$ and consider 
the fundamental domains
of $F$ given by  
$D^+=[\beta^{-1}\rho, \rho]$ and
$D^-=[1-\rho,\lambda(1-\rho)]$
contained in the neighborhoods of $0$ and $1$
where $F$ is affine.
We can choose $\rho$ such that 
$F^N(D^+)=D^-$ for some
$N\in \NN$ (for notational
simplicity  let us put $N=1$).
The map $F^N$
is the {\em transition  from  
$0$ to $1$
\footnote{This transition plays a key role for determining the 
dynamics after unfolding the cycle and it  is determined 
by the {\em  Mather invariant\/} of $F$, \cite{Ma}: 
in a neighborhood of $0$, the map $F$ is the exponential
of the vector field
$X(y)=(\log\beta)\, y \, \frac\partial{\partial y}$, 
whose flow is $y\mapsto \beta^t\, y$. Similarly,
in a neighborhood of $1$, $F$ is the time-one of 
$Y(y)=(\log\lambda)\, (y-1)\, \frac\partial{\partial y}$. 
Consider, for  $y$ close to $0^+$,  
an $n$ large such that $F^n(x)$ is close to $1^-$, 
and write $DF^n(y)(X(y))=\mu(y)\, Y(F^n(y))$. 
Using the local $F$-invariance of $X$ and $Y$ (near $0$ and $1$), 
one has that
$\mu(x)$ does not depend on $n$ and that 
$\mu(x)=\mu(F(x))$. The function $\mu$ is the {\em Mather invariant\/}
 of $F$,
which describes its distortion. For instance, 
if $\mu$ is identically $1$, then $F$ is exactly the exponential of a
vector field.}\/}.   
Suppose that that the eigenvalues $\lambda$ and $\beta$ and
the map $F$  verify the following conditions:

\begin{description}
\item{{\bf (T1)}} 
$\displaystyle{
F^\prime(x) \geq 
\frac{1}{2} \, \frac{1 - \lambda}
{1 - {\beta}^{-1}}}, 
\mbox{ for all } x \in D^+$,
\item{{\bf (T2)}} 
$ (1 - \lambda) < {\beta}^{-1},$ and
\item{{\bf (T3)}} 
$ \displaystyle{ \frac{(1 - \lambda)\,\lambda}{2\,(1 - \beta^{-1})\,\beta}
= \ell  > 1.}$
\end{description}
To get condition (T1) it is enough to consider 
$F$ with small distortion.
To have conditions (T2) and (T3) it is enough to take
$\beta$ close enough to $1^+$.
In fact, later we will consider the case where
$Q$ is a saddle-node (saddle-node heterodimensional cycles),  
see Section~\ref{s.selano}.
Our first result is:

\begin{theor}
\label{t.lat}
For every  $t>0$ sufficiently small,  
$H(P,f_t)\subset H(Q,f_t)$ and $\La_t\subset H(Q,f_t)$.
\end{theor}

This theorem was stated in \cite{Detds}. Here we give a more conceptual
prove of it, which enables us to introduce some technical
tools to be used systematically later on. 
First, in Section~\ref{s.dinamica1dim} we will introduce the system of iterated
functions associated to the cycle (this approach is motivated by
\cite{DRfund}).
In Section~\ref{s.pontosnaoerrantes}, 
we deduce the theorem
from the results in Section~\ref{s.dinamica1dim}.

\section{Expanding one-dimensional dynamics associated to the cycle}
\label{s.dinamica1dim}

For each small  $t>0$,  consider the  
{\em scaled fundamental  domains\/} $D_t^\pm$ defined as follows: 
let
$D^-_t = [1-t,\lambda (1-t)]$ 
and define
$k_t$ as the smaller $k \in \NN$ with
$F^{-k}(D^-_t) \subset [0,t]$. We define
$$
D_t^+ = [a_t,b_t] = 
[\beta^{-1} (b_t),b_t]=
 F^{-k_t}(D^-_t),
\quad
\mbox{where}
\quad
\beta^{-2} t < a_t<b_t\le t.
$$

We next define an expanding  map
$R_t$, 
$R_t\colon D^+_t\to D_t^+$, with discontinuities that will   
describe
 the dynamics 
(in the central 
direction) of the return map of  $f_t$ defined on 
$[-1,1]\times D_t^+ \times [-1,1]$. For each small $t>0$
define the {\em transition map $T_t$ from  $D_t^+$ to $D_t^-$\/} by
$$
T_t \colon  D_t^+ \rightarrow D_t^-,
\quad
x\mapsto T_t(x)=F^{k_t}(x).
$$
\begin{lemm}
\label{l.Tt}
The map 
$T_t$
verifies $T_t^\prime(x)>\ell>1$ for all $x\in D_t^+$, 
where
$\ell$ is as in condition {\bf (T3)}.
\end{lemm}

\begin{demo}
Given $x\in D_t^+$ let 
$k_t=n_t(x)+1+m_t(x)$, where
$F^{n_t(x)}(x)\in D^+$ and 
$F^{i}(x)\not \in D^+$ for all $1 \leq i < n_t(x)$.
We claim that
\begin{equation}
\label{e.blnt} 
\frac {1}{\beta \, t} \leq {\beta}^{n_t(x)} \leq \frac{{\beta}^2}{t}
\quad
\mbox{and}
\quad
\lambda\,  t \le {\lambda}^{m_t(x)} \le \frac{t}{\lambda}.
\end{equation}
For first pair of inequalities just note 
$x\in D_t^+ \subset (\beta^{-2} \, t, t]$ 
and
$\beta^{n_t(x)}\, x \in [\beta^{-1},1]$.
The other ones follow
analogously.

Since $T_t(x) = \la^{m_t(x)} \, F (\be^{n_t(x)}\,x)$,
hypotheses (T1) and (T3) and the estimates in 
(\ref{e.blnt}) give
$$
|T_t^\prime (x)| = \la^{m_t(x)}\, |F^\prime (x)| \, {\be}^{n_t(x)} \geq
(\lambda \, t) \,
\left( \frac{1}{2} \,
\frac{1 - \lambda}{1 - {\beta}^{-1}}\right) \,
\left( \frac{1}{\beta \,t} \right) \, 
= \frac{1}{2} \, \frac{\la\, (1 - \la)}{\beta\, (1 - \beta^{-1})} = \ell > 1,
$$
as we claimed.
\end{demo}

\medskip
Since $T_t(x) \in D_t^-=[1-t,1-\lambda t]$ for all $x \in D_t^+$,
we can define the map $G_t$ by
$$
G_t \colon D_t^+ \rightarrow [0,t],
\quad
x \mapsto G_t(x) = T_t(x) + (t-1).
$$

\begin{rema}
\label{r.Gt}
The map $G_t$ is monotone increasing and
$G_t(D_t^+)=[0,t(1-\lambda)]$.
\end{rema}

\begin{clai}
\label{c.ai>0}
Let $(a_t,b_t] = \tilde{D}_t^+\subset D_t^+$.
Given  $x\in \tilde{D}_t^+$ let $i(x)\in \ZZ$ 
be the minimum $i$ with
${\beta}^i(G_t(x)) \in D^+_t$. 
Then there is $i_0>0$ (maximum with such property)
such that $i(x)\ge i_0$ for all $x \in\tilde D_t^+$. 
\end{clai}

\begin{demo}
Observe first that, since
$\tilde{D}_t^+ \subset [{\beta}^{-2}t, t]$,  
$b_t\in (\beta^{-1}\, t,t]$.
On the other hand, from Remark~\ref{r.Gt}
and $(1 - \la) < \beta^{-1}$ (condition (T2)), 
$G_t(\tilde{D}^+_t) = (0, t\,(1 - \lambda)]
\subset
(0, \beta^{-1}\, t)$. 
Thus 
the right extreme of 
$G_t(\tilde{D}_t^+)$ is less than  the
left extreme of $D_t^+$, hence $i(x)>0$ for all $x\in \tilde D_t^+$,
ending the proof of the claim.
\end{demo}

Finally,
the return map $R_t$ is defined by
$$
R_t\colon \tilde{D}^+_t \rightarrow D^+_t, 
\quad
R_t(x)={\beta}^{i(x)}(G_t(x))=\beta^{i(x)}\, (T_t(x)+(t-1)).
$$
Next we study the dynamics of $R_t$: the
map $R_t$ 
is uniformly expanding and has
(infinitely many) 
discontinuities 
where the lateral derivatives are well defined. 
These discontinuities will play a key role in our constructions.   
The definition of 
$i_0\in \NN$ in Claim~\ref{c.ai>0}
implies that
$\beta^{-i_0}(a_t) \in G_t(D_t^+)$.  
For each $i\ge i_0$ define
$d_i\in \tilde{D}_t^+$ by  $G_t(d_{i}) = \beta^{-i}(a_t)$. 
By construction, the sequence $(d_i)_{i\ge i_0}$
corresponds to the discontinuities of $R_t$ and
verifies the following:
\begin{itemize}
\item
$d_{i+1}< d_i$ and
$d_i\to a_t$, 
\item
Let $[d_{i+1}, d_i]=I_i$, $i > i_0$, and
$I_{i_0} = [d_{i_0},b_t]$.
The map $R_t$ is continuous and 
strictly increasing in the interior of each interval $I_i$.
We continuously extend $R_t$ to the whole $I_i$,
obtaining a bi-valuated return map $R_t$ with
$R_t(d_i)=\{a_t, F(a_t)=b_t\}$ for all $i>i_0$.
In particular, the restriction of $R_t$ to any $I_i$, $i>i_0$, 
is onto.
We let $R_t(b_t)=c_t\le b_t$.
\end{itemize}
The main properties of $R_t$ are summarized in the next lemma.

\begin{lemm}
\label{l.Rtexpansora}
The restriction of $R_t$ to each interval
$I_i, i > i_0$, is onto and 
$R_t^\prime(x)>\ell >1$ for all $x\in (a_t,b_t]$
(if $x=d_i$ this means that the lateral derivatives of
$R_t$ at $x$ are  greater than $\ell$).
Moreover,
$G_t(R_t(d_i)) = 0$
for all $i \geq i_0$.
\end{lemm}
The  expansiveness
of $R_t$ follows 
from Lemma \ref{l.Tt} and Claim~\ref{c.ai>0}. Condition 
$G_t(R_t(d_i))=0$ follows from
$R_t(d_i)=a_t$ and
$G_t(a_t)=0$.

\begin{lemm}
\label{l.descont}
Consider a small $t>0$ 
and   an open subinterval $J$  of $\tilde{D}^+_t$. Then there is  
a first $k \in \NN \cup \{0\}$ 
such that $R_t^k(J)$  contains a discontinuity of $R_t$.
In particular, there is
$x\in J$ 
such that
$G_t(R_t^k(x))=0.$
\end{lemm}

\begin{demo}
If the interval  $J$ contains a discontinuity we are done.
Otherwise, let $i>0$ be such that the intervals $J$, $R_t(J), 
\dots, R_t^i(J)$ do not contain
discontinuities. Thus,
for each $k\in \{0,\dots, i\}$, there is $i_k\ge i_0$ such that
$R_t^k(J)\subset I_{i_k}$.
Lemma~\ref{l.Rtexpansora}
implies that
$|R_t^k(J)|\ge \ell^k |J|$, $\ell>1$,
for all  $k\in \{0,\dots, i\}$.
Since the size of the intervals $I_i$ is upper bounded,
this inequality implies that
there is a first $m\in \NN$
such that
$R_t^m(J)$ 
is not contained in any $I_i$, thus
it intersects the set of  discontinuities of $R_t$.
\end{demo}

\section{The maximal invariant set: Proof of Theorem~\ref{t.lat}}
\label{s.pontosnaoerrantes}

Next proposition is the main technical result of this section. 
Heuristically, it 
 means that the one-dimensional stable 
manifold of $Q$ topologically behaves as a 
two-dimensional manifold.

\begin{prop}
\label{p.chi}
For every 
small $t>0$ and every
two-disk $\chi$ with 
$W^s(P,f_t)\pitchfork \chi\ne \emptyset$,
$W^s(Q,f_t)\pitchfork \chi\ne \emptyset$.
In particular, $W^s(P,f_t)$ is contained
in the closure
of $W^s(Q,f_t)$.         
\end{prop}

\paragraph{Proof of the  inclusion $H(P,f_t)\subset H(Q,f_t)$ 
in Theorem~\ref{t.lat}.}
This inclusion follows from 
Proposition~\ref{p.chi}. 
By the definition of $H(P,f_t)$, it suffices 
to see that any  
$x\in W^s(P,f_t)\pitchfork W^u(P,f_t)$
is accumulated by homoclinic points of $Q$.
The geometric configuration of the cycle implies that
$W^u(P,f_t)\subset
\mbox{closure\,} (W^u(Q,f_t))$. Thus,
given any $x\in H(P,f_t)$
and any $n>0$, there is a disk
$\Delta_n$ simultaneously contained in $W^u(Q,f_t)$
and in the ball of  radius $1/n$ centered at  
$x$ whose   
interior meets transversely  $W^s(P,f_t)$.
By Proposition~\ref{p.chi},  
$\De_n\pitchfork W^s(Q,f_t)\ne \emptyset$. Thus there is $y_n\in
\De_n\cap H(Q,f_t)$. 
By construction, $y_n\to x$,
proving the inclusion
$H(P,f_t)\subset H(Q,f_t)$.

\medskip

\subsection{Homoclinic classes - Proof of Proposition~\ref{p.chi}}
\label{ss.homoclinicclasses}

We now go into the details of the
proof of Proposition~\ref{p.chi}. We
first introduce some definitions.

\begin{itemize}
\item
A set
$\Delta \subset [-1,1]\times [-1,2]\times [-1,1]$ 
is a {\em vertical strip\/}
if 
$\Delta = \{x_1\}\times[l_1,l_2]\times[r_1,r_2]$,
where $l_1<l_2$ and  $r_1<0<r_2$.
The segment 
$\{x_1\}\times [l_1,l_2]\times \{0\}$ is the {\em basis\/} of $\De$. The
{\em width\/} and the 
{\em height\/} of $\Delta$ are
 $w(\Delta)=(l_2-l_1)$ and 
$h(\Delta)=(r_2-r_1)$. 
The strip $\De$ is {\em complete\/} if
$r_2=1$ and $r_1=-1$, 
{\em well located\/} if
$[l_1,l_2]$ is contained in the interior of $D_t^+$,
and {\em perfect\/}
if it simultaneously is complete and well located.
\item
A subset 
$J\subset  [-1,1]\times [-1,2]\times [-1,1]$ 
is a {\em vertical segment\/} if  
$J = \{x_1\}\times\{l_1\}\times[r_1,r_2]$, where
$r_1<0<r_2$.
The  {\rm height\/} of $J$ is   $h(J)=(r_2-r_1)$. 
As above, the  segment $J$ 
is {\em complete\/} if  
$r_1=-1$ and $r_2=1$, {\em well located\/}
if $l_1$
 is in the interior of
$D_t^+$, and {\em perfect\/} if it simultaneously is complete and well
located.
\item
A vertical segment $J$ 
(resp. strip $\De$) 
is {\em at the right
of $Q$\/} if
$l_1\in (0,1]$. 
\end{itemize}

Given an interval $\alpha \subset [-1,2]$,
denote by $\De(\{x\}\times \alpha \times \{0\})$
(resp. $J (x,y,0)$) the {\em unique\/}
complete vertical strip (resp. segment)
with basis $\{x\}\times \alpha\times \{0\}$
(resp $(x,y,0)$).
The following algorithm associates to perfect segments and strips
their successors:

\begin{algo}
\label{a.algoritmo}
Let 
$\De=\De(\{x\}\times \alpha \times \{0\})$
be a perfect strip and define 
$\cG_t(\De)$  as
the  perfect strip
such that:
\begin{itemize}
\item
the basis of
$\cG_t(\De)$ 
is of the form $(\{x^\prime\}\times 
G_t(\alpha) \times \{0\})$, where $x^\prime=\la_s^{k_t} \, x-1/2$,
\item
$\cG_t(\De)$ is contained in $f_t^{k_t}(\De)$
(where $F^{k_t}(D_t^+)= D_t^-$).
\end{itemize}
Suppose now that $\alpha$ does not contain discontinuities, i.e.
$\alpha\subset (d_i,d_{i+1})$ for some $i$. 
Define $\cR_t(\De)$  as
the  perfect strip
such that:
\begin{itemize}
\item
the basis of
$\cR_t(\De)$ 
is of the form $(\{x^\prime\}\times 
R_t(\alpha) \times \{0\})$,  where
$x^\prime =\la_s^i(\la_s^{k_t} \, x-1/2)$,
\item
$\cR_t(\De)$ is contained in $f_t^{k_t+i}(\De)$.
\end{itemize}
Similarly, to a perfect segment $J =J (x,y,0)$ 
we associate
 perfect segments $\cG_t(J)$ and $\cR_t(J)$
(provided $y\ne d_i$ for all $i$).
\end{algo}

The strips $\cG_t(\De)$ and $\cR_t(\De)$ in Algorithm~\ref{a.algoritmo}
are obtained as follows.
Given a set $A$ and a point $x\in A$,  denote by
$C(x,A)$ the connected component of $A$ containing $x$. 
Consider a small tubular neighborhood $V$ of
$f^{k_0}_0(\{0,1\}\times [-1,1])\subset W^u(P,f_0)$,
$k_0$ as in 
the definition of the cycle in Section~\ref{s.heterodimensional},
then 
$$
\cG_t(\De)=
C(f^{k_t}_t(x,y,0), f_t^{k_t}(\De) \cap V)
\cap [-1,1]^3),
\qquad
\cR_t(\De)=
C(f^{i}_t(x^\prime,y^\prime,0), f_t^{k_t+i}(\De) \cap V)
\cap [-1,1]^3),
$$
where
$(x,y,0)$ is any point in the basis of $\De$, 
$x^\prime=(\la_s^{k_t} \, x-1/2)$, and  $y^\prime\in G_t(\alpha)$.
The construction for the successors of the segments is analogous.

\begin{lemm}
\label{l.Wuptsegtovert}
The manifold $W^u(P,f_t)$ contains a perfect 
segment for all small  $t>0$.
\end{lemm}

\begin{demo}
Consider the transverse homoclinic point 
$x_t=(-1/2,t,0)$ of
$P$ (for $f_t$).
Recall that
$t\ge b_t$ and $\beta^{-1} \,t\in D_t^+=[a_t,b_t]$.
Let us assume that 
$t> b_t$,
and thus $\beta^{-1} t\in (a_t,b_t)$
(the case $t=b_t$ follows similarly, so it will be omitted).
Consider the complete vertical segment $\cR_t(J)$, where 
$J=
J(-(\la_s^{-1}/2), \be^{-1}t,0)\subset W^u(P,f_t)$.
If $R_t(\be^{-1} t)$ belongs to the interior of $D_t^+$,
then $\cR_t(J)\subset W^u(P,f_t)$ is
the announced segment.
Otherwise,
$R_t(\be^{-1} t)=b_t$ and there is
a homoclinic point of
$P$ of the form $(x^\prime, b_t,0)$.
Using the  $\la$-lemma and the 
product structure of the cycle, 
one gets  
homoclinic points  $(x_n, y_n,0)$ of $P$ 
and  complete
segments $J_n=J(x_n, y_n,0)\subset W^u(P,f_t)$
such that 
$x_n\to x^\prime$, $y_n\to b_t$, and 
$y_n$ is increasing.
 Thus $y_n$ belongs to the interior of $D_t^+$ for every
big $n$ and $J_n\subset W^u(P,f_t)$ is perfect.
\end{demo}

For clearness we first prove
Proposition~\ref{p.chi} 
in the following special case:

\begin{prop}
\label{p.chip}
Let $\chi \subset [-1,1]\times [-1,2]\times [-1,1]$ be a set of the form 
$\{x\} \times A$, where $A$ is a disk 
of $\RR^2$ whose interior contains a point of the form $(y,0)$ with
$y\in (0,2)$.
Then $\chi$  intersects transversely $W^s(Q,f_t)$.
\end{prop}

We claim that is enough to prove the 
Proposition~\ref{p.chip} for perfect strips:

\begin{lemm}
\label{l.faixavertical}
Let $\De$ be a perfect
strip. 
Then there is $k \in \NN$
such that
$f_t^{k}(\De) \pitchfork W^s(Q,f_t) \neq \emptyset.$
\end{lemm}

\begin{demo}
Suppose that $\De=\De(\{x_0\}\times \alpha\times \{0\})$, 
$\alpha$ in the interior of $D^+_t$. 
By Lemma~\ref{l.descont}, 
there exist $y_0$ in the interior of $\alpha$ 
and  $k\in \NN$
such that $G_t(R^k_t(y_0))=0$.
Thus the vertical strip
$\cG_t(\cR^k_t(\De))$ (contained in the forward orbit of
$\De$) intersects transversely  $[-1,1]\times \{(0,0)\}
\subset W^s(Q,f_t)$.
\end{demo}

\medskip
\noindent
{\bf Proof of 
Proposition~\ref{p.chip}:}
By Lemma~\ref{l.Wuptsegtovert},
$W^u(P,f_t)$
contains a perfect vertical segment $J$.
By definition, $\chi$ meets transversely 
$W^s(P,f_t)$, which implies, by the
$\lambda$-lemma, that the
forward orbit of 
$\chi$ contains a sequence of complete  strips $\chi_n$
accumulating to $J$. 
Thus $\chi_n$
contains a perfect  strip for all  $n$ large. 
Lemma~\ref{l.faixavertical}
now implies that $\chi_n$  (and thus $\chi$) 
transversely meets $W^s(Q,f_t)$, ending the proof of the proposition.
\hfill$\Box$

\medskip
\noindent
{\bf Proof of 
Proposition~\ref{p.chi}:}
We can  assume that
$\chi$ is  transverse to
$W^s_{loc}(P,f_t)$ and contained in 
$[-1,1]\times [-1,2]\times [-1,1]$.
If  $\chi$ contains a 
subset of the form $\{x\}\times A$, where  
$A$ is an open subset of $\RR^2$ containing a point  $(0,y)$
with $y\in (0,2)$,    
Proposition~\ref{p.chip} implies the result.
For the general case, consider
 a point $(x_0,y_0,0)$, $y_0\in (0,2)$, in the interior of
$\chi\pitchfork W^s_{loc}(P, f_t)$ and for every big $n$
the  vertical strip  
$$
\Sigma_n=\{x_0\}\times
[y_0-1/n, y_0+1/n]\times
[-1/n,1/n].
$$
The strips $\Si_n$ 
verify the hypotheses of 
Proposition~\ref{p.chip},
hence there is $(x_0,y_n,z_n)\in \Sigma_n \pitchfork W^s(Q,f_t)$
such that 
$H_n=
[-1,1]\times \{(y_n,z_n)\}
\subset W^s(Q,f_t)$.
Since 
$(x_0,y_n,z_n)\to (x_0,y_0,0)$,
it is immediate that 
$H_n$ meets transversely  $\chi$ for all large $n$, ending the proof of the
proposition.
\hfill$\Box$

\subsection{The maximal $f_t$-invariant set}
\label{ss.maximalinvariant}

To prove the second 
part of Theorem~\ref{t.lat} 
($\La_t\subset H(Q,f_t)$)
let
$V_0$ be the connected component of 
the neighborhood of the cycle $W$ containing the heteroclinic
point $(-1/2,0,0)$.
There are two types 
of points of $\La_t$: 
(a) those points whose orbit
does not meet $V_0$ (i.e., the set 
$\{0\}\times [0,1]\times \{0\}$) 
and (b) those having an iterate in $V_0$.

We claim that
every point of type (a) belongs to $H(Q,f_t)$: given 
any $(0,x,0)$, $x\in (0,1)$, consider the disk 
$\De_n
=\{0\}\times [x-1/n,x+1/n]\times [-1/n,+1/n]
\subset W^u(Q,f_t)$ satisfying the hypothesis of  Proposition~\ref{p.chi}. 
Hence $\De_n\pitchfork W(Q,f_t)\ne \emptyset$ and thus
$\De_n \cap H(Q,f_t)\ne \emptyset$. Since this holds for all $n\in \NN$,
$(0,x,0)\in H(Q,f_t)$. 

For points $w\in \La_t$ of type (b),
after replacing $w$ by
some iterate of it,  we can assume that $w\in V_0$.
Consider the sequence
$(n_i(w))_{i \in \cI(w)}$ associated to $w$,  where 
$\cI(w) \subset \ZZ$ is an interval in $\ZZ$,
inductively defined as follows:
let $n_0(w)=0$ and,  assuming defined $n_j(w)$, $j\ge 0$,
we define
$n_{j+1}(w)$ as the first integer $k>n_j(w)$ such that
$f_t^k(w) \in V_0$.
We argue analogously for negative $j$, assuming defined $n_j(x)$, 
$j\le 0$,
$n_{j-1}(w)$ is the first negative
integer  $k<n_j(w)$ with
$f_t^k(w) \in V_0$. 

Define 
$\cI^+_t(b)$ 
as the subset of 
$\La_t\cap V_0$
of points $w$ such that
$\cI(w)$ is  upper bounded. 
The subset
$\cI^-(b)$  
is defined similarly.
We let 
$\cI^\pm (b)=\cI^+(b)\cap 
\cI^-(b)$.
We borrow from
\cite[Lemma 4.1]{DRnon} the following lemma whose proof is 
straightforward:

\begin{lemm} 
\label{l.classificacao}
For every small $t>0$,   
$\cI_t^+(b)\subset W^s(P,f_t)\cup W^s(Q,f_t)$ and
$\cI_t^-(b)\subset W^u(P,f_t)\cup W^u(Q,f_t)$. 
\end{lemm}

Next result 
immediately follows by  
observing that 
$f_t$ (resp. $f_t^{-1}$)
exponentially expands the vertical (resp. horizontal) segments:

\begin{rema}
\label{r.expansaonegativa}
Let $w=(x,y,z)\in \cI_t^+(\infty)$ (resp.
$w\in \cI_t^-(\infty)$).
Then 
$\{(x,y)\}\times [z-\ve,z+\ve]
\pitchfork W^s(P,f_t)\ne \emptyset$
(resp. 
 $[x-\ve,x+\ve]\times \{(y,z)\}
\pitchfork W^u(Q,f_t)\ne \emptyset$)
for every $\varepsilon>0$. 
\end{rema}

To prove the theorem we consider the following four cases.  

\medskip

\noindent
{\bf Case (i):
$w=(x,y,z)\in \cI^-_t(b) \setminus \cI^+_t(b)$.} 

By Remark~\ref{r.expansaonegativa},
there is a sequence 
$w_n=(x,y,z_n)\in W^s(P,f_t)$ with $w_n\to w$.
We claim that $w_n\in H(Q,f_t)$ for all large $n$.
Thus $w\in H(Q,f_t)$.
To prove the claim, note that
the distances between the backward iterates of
$w_n$ and $w$ exponentially decrease,
we get that $w_n\in \La_t$. 
Moreover, since $w\in \cI^-_t(b)$,
Lemma~\ref{l.classificacao}
implies that
 $w, w_n\in W^u(P,f_t)\cup W^u(Q,f_t)$.
If $w_n\in W^u(P,f_t)$, then $w_n\in H(P,f_t)\subset
H(Q,f_t)$ (recall the first part of Theorem~\ref{t.lat}
proved above) and we are done.
Otherwise,  
$w_n\in W^u(Q,f_t)$ and 
for each  $k$ large, there is a small vertical
strip $\De_k$ of diameter less than $1/k$, 
whose interior is contained in 
$W^u(Q,f_t)$ and contains $w_n$.
Since $w_n\in W^s(P,f_t)$,
$\De_k\pitchfork W^s(P,f_t)$.
Thus,
by Proposition~\ref{p.chi}, 
$W^s(Q,f_t)$ intersects transverselly  the interior of
$\De_k$.
Hence,
since the interior of  $\De_k$ is contained in $W^u(Q,f_t)$, 
$\De_k$ contains a homoclinic
point $y_k$ of $Q$. From
$\mbox{diam}(\De_k)\to 0$, we get $y_k\to w_n$, which implies
 $w_n\in H(Q,f_t)$.

\medskip
\noindent
{\bf Case (ii):
$w=(x,y,z)\not \in \cI_t^+(b) \cup \cI_t^-(b) $.}  

We claim that $w$ is accumulated by points
$w_n\in \cI_t^+(\infty)\cap \cI^-_t(b)$, and 
the result follows from Case (i).
To prove the claim observe that,
by Remark~\ref{r.expansaonegativa},
there is a sequence 
$w_n=(x_n,y,z)\in W^u(Q,f_t)$ with $w_n\to w$.
Since the distances between the forward iterates of
$w_n$ and $w$ exponentially decrease,
$w_n\in \La_t$. 
This also implies that
$w_n\in \cI^+_t(\infty)$.
Finally, 
$w_n\in W^u(Q,f_t)$ implies 
$w_n\in \cI^-_t(b)$, ending the proof of the claim.

\medskip
\noindent
{\bf Case (iii):
$w=(x,y,z)\in \cI_t^+(b)\setminus \cI^-_t(b)$.}  

By Lemma~\ref{l.classificacao},
 $w \in W^s(P,f_t)\cup W^s(Q,f_t)$ and, by replacing $w$ by a forward iterate,
 we can assume that $w=(x,y,0)$, $y\ge 0$. 
Remark~\ref{r.expansaonegativa} gives a sequence 
$w_n=(x_n,y,0)\in W^u(Q,f_t)$ with  $w_n\to w$.
For each $n$, there is a vertical disk $\De_n\subset  W^u(Q,f_t)$ 
centered at  $w_n$, of 
diameter less than  $1/n$. Clearly, $\De_n$ intersects transversely 
$W^s(P,f_t)$. Thus, by
Proposition~\ref{p.chi},   
$\De_n \pitchfork W^s(Q,f_t) \neq \emptyset$. As in the previous cases, this
implies that
$\De_n \cap H(Q,f_t) \neq \emptyset$ for all  $n$ large,
thus $w\in H(P,f_t)$.

\medskip
\noindent
{\bf Case (iv):
$w\in \cI^\pm_t(b)$.}

By  Lemma~\ref{l.classificacao},
there are four possibilities:
(1) $w \in W^s(Q,f_t) \cap W^u(Q,f_t)$,
(2) $w \in W^s(P,f_t) \cap W^u(P,f_t)$,
(3) $w \in W^s(P,f_t) \cap W^u(Q,f_t)$,
and (4)
 $w \in W^u(P,f_t) \cap W^s(Q,f_t)$.
Recall that the intersections above are transverse or quasi-transverse, 
depending on the case. Hence, 
in case (1), $w\in H(Q,f_t)$ and, in
case (2), $w\in H(P,f_t)\subset H(Q,f_t)$. 
In case (3),
the  same proof of $\{0\}\times [0,1]\times \{0\}\subset H(Q,f_t)$ 
implies that
$w\in H(Q,f_t)$: 
just observe that for every disk $\De\subset W^u(Q,f_t)$ containing $w$,
$W^s(P,f_t)\pitchfork \De\ne 
\emptyset$, thus $\De\cap H(Q,f_t)\ne \emptyset$.
It still remains the case
$w \in W^s(Q,f_t) \cap W^u(P,f_t)$.
By replacing $w$ by a forward iterate, we can assume
that $w=(x,0,0)$, $x\in [-1,1]$, and the following lemma 
and Cases (i) and (ii) easily imply 
case (4):

\begin{lemm}
\label{l.zninfty}
Let $w=(x,0,0) \in V_0 \cap (W^s(Q,f_t) \cap W^u(P,f_t))$.
Then there is a  sequence
$w_n\to w$ with 
$w_n\in \cI_t^+(\infty)$.
\end{lemm}

\begin{demo}
For each $n\in \NN$,
consider the rectangle
$R_n(x) = \{x\}\times [0, 1/n]\times [-1/n,1/n]$.

\begin{clai}
\label{c.kn}
There exists $\kappa_n$ in  $R_n(x)\in \La_t^+$ 
whose forward orbit returns to  $V_0$ infinitely many times.
\end{clai}

Assuming this claim, we now finish the proof of the lemma:
similarly as in Remark~\ref{r.expansaonegativa},
but now considering points in $\La_t^+$, 
we have that the point  $\kappa_n=(x_n,y_n,z_n)$ is 
accumulated by points $\kappa_n^m\in W^u(Q,f_t)$
of the form $(x_n^m, y_n, z_n)$. Since the distances between the forward 
iterates of $\kappa_n$ and $\kappa_n^m$ decrease, 
the forward orbit of $\kappa_n^m$
is contained in $W$ and returns infinitely many times to $V_0$. 
On the other hand,
since $\kappa_n^m\in W^u(Q,f_t)$, its backward orbits
also is in $W$. Thus the whole orbit of $\kappa_n^m$ is in $W$, so
$\kappa_n^m\in \La_t$ and $\kappa_n^m\in \cI_t^+(\infty)$.
By Cases (i) and (ii) above, $\kappa_n^m\in H(Q,f_t)$, thus
$\kappa_n\in H(Q,f_t)$.
To prove the claim, we  need the following fact:

\begin{fact}
\label{l.retangulo}
Let $R=R_n(x)$, 
$1/n<t$. Then 
there is
$i=i(R) \in \NN$ 
such that, for every 
$j \geq i$,
$f_t^j(R)$ contains a  rectangle
$\Ga(R,j) $ of the form
$\{a\}\times [0,1/n]\times [-1/n,1/n]$. 
 \end{fact}

\begin{demo} 
Let
$N_t$ be the smaller $i \in \NN$ such that
$F^{i}(1/n)\in (1-t+1/n,1)$
and write
$$
e=(1-t+1/n+g)=F^{N_t}(1/n), 
\quad
g\in (0,t-1/n).
$$
The definition of the unfolding of the cycle implies that,
for each $j\ge 0$,
$f_t^{N_t+j}(R_n(x))$ contains a rectangle of the form
$$
\{\la_s^{N_t+j}(x)-1/2\}
\times 
[0,t+\la^j(g+1/n-t)]
\times [-1,1]
\supset
\{\la_s^{N_t+j}(x)-1/2\}
\times 
[0,1/n]
\times [-1,1],
$$
where the inclusion follows from
$
t+\la^j(g+1/n -t) \ge
t +\la^j(1/n -t)\ge 1/n.
$
This finishes the proof of the fact.
\end{demo}

To prove Claim~\ref{c.kn},
consider   $R_n(x) = R(0)$ and, using Fact~\ref{l.retangulo}, let 
$R(1)=\Ga (R(0), i(R(0))$.
Write $i_0=i(R_0)$ and
$R^1 =  f_t^{-i_0}({R}(1)) \subset R(0)$.
Assume inductively defined
numbers $i_{k-1}$ and  
rectangles $R(k)$ and $R^k$ for every $k\in \{0,\dots,j\}$
as follows:
\begin{itemize}
\item 
$R(k)=\Ga(R(k-1), i(R(k-1)))$
and $i_{k-1}=i(R(k-1))$,
in particular,
$R(k)$ satisfies the hypotheses of Fact~\ref{l.retangulo}, 
\item
$R^{k}\subset R^{k-1}\subset \cdots \subset R^1\subset R(0)=R_n(x)$
and $R^k=f_t^{-i_0-\cdots -i_{k-1}}(R(k))$.
\end{itemize}
We define
$i_k=i(R(k))$,
$R(k+1)=\Ga(R(k), i_k)$ and
$R^{k+1}=f_t^{-i_0-\cdots -i_{k}}(R(k+1))$,
completing the inductive process.
Now it suffices to take any point in the non-empty intersection
$\cap_{k\in \NN} R^k$.
\end{demo}

The proof of Theorem~\ref{t.lat} is now complete.

\section{Saddle-node heterodimensional cycles}
\label{s.selano}

In this section, we consider 
{\em saddle-node heterodimensional cycles.\/}
For that, in the definition of the heterodimensional cycle in 
Section~\ref{s.heterodimensional}, we replace the function $F$ (defining the
central dynamics) by a one parameter family of 
maps 
$\Phi_s\colon [-1,2]\to \RR$ such that: 
\begin{itemize}
\item
For every $s$, 
the point $1$ is an attracting hyperbolic point of $\Phi_s$ 
and $\Phi_s$ is linearizable in a neighborhood
of $1$ (independent of $s$). We denote by $0<\la<1$ the eigenvalue
of $\Phi_s$ at $1$.
\item
Locally in $0$, the map $\Phi_s$ is of the form
$\Phi_s(x)= x+x^2 -s$.
Thus, for $s>0$, $\Phi_s$ has
two hyperbolic fixed points
$\pm \sqrt{s}$ (an attractor and a repellor)
collapsing at $s=0$. Moreover,  for every $s<0$,
$\Phi_s$ has no fixed points close to $0$. 
\item
Every $\Phi_s$ is strictly increasing and
has no fixed points different from
$1$ and $\pm\sqrt{s}$.
\end{itemize}


We now define, as in Section~\ref{s.heterodimensional}, a two parameter
family of diffeomorphisms
$f_{t,s}$: 
the parameters $t$ and $s$ describing 
the motion of the unstable manifold
of $P$ and the unfolding of the
saddle-node  (i.e., $f_{t,s}(0,y,0)=(0, \Phi_s(y),0)$), 
respectively.
Observe that 
$P=(0,1,0)$, $S^-_s=(0,-\sqrt{s},0)$ and $S^+_s=(0,\sqrt{s},0)$
($s\ge 0$)
are fixed points of
$f_{t,s}$.

We let $f_t=f_{t,0}$.
For the saddle-node  $S=(0,0,0)$ of $f_t$ there are
defined the stable and unstable manifolds 
(denoted $W^s(S,f_t)$ and $W^u(S,f_t)$)
and the strong stable and unstable manifolds
(denoted by $W^{ss}(S,f_t)$ and
$W^{uu}(S,f_t)$).
Observe that
$W^s(S,f_t)$ and $W^u(S,f_t)$
are
two-manifolds with boundary and 
$W^{ss}(S,f_t)$ and
$W^{uu}(S,f_t)$ have 
both dimension one. 
Notice that
$$
\begin{array}{ll}
&\{0\}\times [0,1)\times [-1,1]
\subset W^u(S,f_t),
\quad
[-1,1]\times [-2,0]\times \{0\}
\subset W^s(S,f_t),
\\
&[-1,1]\times \{(0,0)\}
\subset W^{ss}(S,f_t),
\quad
\{(0,0)\} \times [-1,1]
\subset W^{uu}(S,f_t).
\end{array}
$$
Keeping in mind these relations, we have that, 
\begin{itemize}
\item
for all $t$,
$W^u(S,f_t)$
meets transversely  $W^s(P,f_t)$ 
throughout the segment $\{0\}\times (0,1)\times \{0\}$,
\item
for $t=0$,
$W^u(P,f_0)$ meets quasi-transversely  $W^{ss}(S,f_0)$
along the orbit of
$(-1/2,0,0)$, 
\item
for $t>0$,
the point $(-1/2,t,0)$ is a transverse homoclinic point of $P$ and 
$(-1/2,0,0)$ is a point of transverse intersection
between $W^{ss}(S,f_t)$ and $W^u(S,f_t)$.
\end{itemize}
In this case, we say that
the arc 
$f_{t}=f_{t,0}$
has a {\em  saddle-node heterodimensional cycle\/}
associated to $P$ and 
$S$ at $t=0$.
This cycle can be thought as a 
a limit case of the heterodimensional
cycles in Section~\ref{s.heterodimensional}, where the derivative of the point
of index two $Q$ is $1^+$.

The two-fold behavior of the saddle-node $S$, as
a point of index two and one simultaneously, 
leads us to consider, for small positive $t$, 
the {\em lateral homoclinic classes of $S$\/} defined by 
$$
H^+(S,f_t)=\overline{W^u(S,f_t) \pfrk
W^{ss}(S,f_t)}
\quad
\mbox{and}
\quad 
H^-(S,f_t)=\overline{W^s(S,f_t) \pfrk
W^{uu}(S,f_t)}.
$$
As in the case of the usual homoclinic classes, we have that:

\begin{prop}
\label{p.lateral}
For every small $t>0$,
$H^+(S,f_t)$ (resp. 
$H^-(S,f_t)$)
is transitive and the periodic points of index
two (resp. one) form a dense subset of it. 
\end{prop}

Consider a neighborhood $W$ of the saddle-node heterodimensional cycle
defined 
as in Section~\ref{s.pontosnaoerrantes} and denote by $\Ups_t$
 the maximal
invariant set of $f_t$ in $W$.

\begin{theor}
\label{t.latlat}
For every small $t>0$, one has that
$H(P,f_t)\subset H^+(S,f_t)$ and $\Ups_t\subset H(S^+,f_t)$.
\end{theor}

The proof of Theorem~\ref{t.latlat}
follows as the one of
Theorem~\ref{t.lat}, the only difficulty being to
redefine appropriately the one-dimensional dynamics
associated to the cycle
(recall Section~\ref{s.dinamica1dim}). 
This will be 
briefly done in the next section.
To get the inclusion
$H^+(S,f_t)\subset H(P,f_t)$
we need the
following distortion
property for the saddle-node map $\Phi=\Phi_0$.

\begin{description} 
\item
{\bf (SN)}
Let 
$K>0$ be the maximum of 
$|\Phi''(x)|/|\Phi'(x)|$,
$x\in [0,1]$,
then
$
\displaystyle{   
\frac{4 \, e^{K}\, (1-\la)}{\la^6}  < \frac{1}{2}}$,
where
$\la\in (2/3,1)$.
\end{description}

\begin{theor}
\label{t.HssubsetHp}
Under the assumption (SN), 
$H^+(S,f_{t})\subset H(P, f_t)$ holds for all small positive 
$t>0$. 
\end{theor}
To prove this theorem we need new ingredients that will be  introduced in
Section~\ref{ss.newingredients}.

\subsection{One-dimensional dynamics for the  saddle-node cycle}
\label{ss.1dim}
We now
adapt the definitions of  scaled fundamental domains,
transitions and returns for saddle-node cycles.
As in Section~\ref{s.dinamica1dim}, for each $t>0$,  
define the
fundamental domains $D_t^-=[1-t, 1-\la\, t]$ 
and $D_t^+=[a_t,b_t]$, $a_t=\Phi^{-1}(b_t)$, 
where $D_t^+$ is the first backward iterate of
$D_t^-$ by $\Phi$ contained in $[0,t]$.
We have 
$\Phi^{k_t}(D_t^+)=D_t^-$, for some $k_t\in \NN$.
Observe that $|D_t^-| = t\, (1 - \la)$ and,
since $b_t\in (0,t]$,
$|D_t^+| \le t^2$.
For small $t>0$,
define the {\em transition\/}
$\mfT_t$ and the map $\mfG_t$ by
$$
\mfT_t \colon  D_t^+ \rightarrow D_t^-,\,\,\,
x\mapsto \mfT_t(x)=\Phi^{k_t}(x)\quad
\mbox{and}
\quad
\mfG_t \colon D_t^+ \to [0, t(1-\lambda)],
\,\,\,
x\mapsto \mfG_t(x)=\mfT_t(x)+t.
$$

\begin{lemm}
\label{l.mfTtexpansora}
The maps $\mfT_t$ and 
$\mfG_t$ are uniformly expanding for all small $t>0$.
\end{lemm}

\begin{demo}
It suffices to see that
$(\Phi^{k_t})^\prime(z)>1$ for all $z\in D_t^+$.
We use the following  standard
lemma (whose proof is omitted here):

\begin{bdl}
\label{f.distorcao}
Let 
$K>0$ be as above.
Then,
for every pair of points 
$z, y\in D_t^+$  and every small $t>0$, it holds 
$$
e^{-K}
\le
\frac{(\Phi^{k_t})^\prime(z)}
{(\Phi^{k_t})^\prime(y)}
\le 
e^{K}.
$$
\end{bdl}
The lemma now follows by the mean value theorem, taking 
$y$ with
$(\Phi^{k_t})^\prime(y)=|D^-_t|/|D_t^+|\ge (1-\la)/t$.
Thus, if $t$ is small,
$(\Phi^{k_t})^\prime(z)\ge 
(e^{-K}\,(1-\la))/t>1$, for all $z\in D_t^+$. 
\end{demo}

As in
Section~\ref{s.dinamica1dim},
given 
$x\in (a_t, b_t]=\tilde D_t^+$, 
let
$i(x) \in \ZZ$ be the first
$i$ with
$\Phi^i(\mfG_t(x)) \in D^+_t.$

\begin{lemm} 
\label{l.left}
There exists $i_0>0$ such that
$i(x)\ge i_0$ for all $x\in \tilde{D}_t^+$.
\end{lemm}

\begin{demo}
To prove the lemma it is enough to see that
$\mfG_t(D_t^+)\subset (0,a_t)$.
By definition,
$\mfG_t(D_t^+) = [0,(1 - \lambda)\, t]$.
Observe that, if $t$ is small enough,
$$
\Phi^2((1-\la) \, t)= 
\Phi((1-\la)\,t +(1-\la)^2\,t^2) = 
(1-\la)\,t +2\, (1-\la)^2\,t^2 + \mbox{h.o.t}
< t.
$$
Thus,
the right extreme $\Phi^2((1 - \lambda) t)$ of 
$\Phi^2(\mfG_t(D_t^+))$ 
is less than
$t$. In particular, the right extreme of $\mfG_t(D_t^+)$
is less than 
$\Phi^{-2}(t)$, and the 
the lemma follows  from $D_t^+ \subset (\Phi^{-2}(t),t]$.
\end{demo}

The {\em  return map\/} $\mfR_t$  is now defined by
$$
\mfR_t\colon \tilde D^+_t \rightarrow D^+_t,
\quad
\mfR_t(x)=\Phi^{i(x)}(\mfG_t(x))=\Phi^{i(x)}(\mfT_t(x)+t).
$$ 
As in the case of the map $R_t$ in 
Section~\ref{s.dinamica1dim}, 
for each $i\ge i_0$, there is
$\de_i\in \tilde{D}_t^+$ with  $\mfG_t(\de_{i}) = \Phi^{-i}(a_t)$.
The points $\de_i$ are the discontinuities of $\mfR_t$. 
In this way, we get 
a decreasing sequence $(\de_i)_{i\ge i_0}$ with $\de_i\to a_t$,
and intervals
$J_i=[\de_{i+1}, \de_i]$, $i > i_0$, and
$J_{i_0} = [\de_{i_0},b_t]$ such that
$\mfR_t$ is continuous and 
increasing in the interior of each $J_i$.
Extending $\mfR_t$ continuously to the whole $J_i$ we
get a bi-valuated map with
$\mfR_t(\de_i)=\{a_t,b_t\}$ for all $i>i_0$.

\begin{lemm}
\label{l.mfRtexpansora}
The restriction of $\mfR_t$ to each interval
$J_i, i > i_0$, is onto. Moreover, there is $\ell>1$ such that 
$\mfR_t^\prime(x)>\ell>1$ for all $x\in (a_t,b_t]$
(if $x=\de_i$ this means that the lateral derivatives of
$\mfR_t$ at $x$ are  greater than $\ell$).
Finally, $\mfG_t(\mfR_t(\de_i)) = 0$
for all $i \geq i_0$.
\end{lemm}

\begin{demo}
The lemma follows as Lemma~\ref{l.mfRtexpansora}
observing that $i_0>0$ (Lemma~\ref{l.left}),
$\mfG_t$ is expanding (Lemma~\ref{l.mfTtexpansora}), 
and that the derivative of
$\Phi$ in $(0,t]$ is bigger than one. 
\end{demo}

Arguing as in 
Section~\ref{s.dinamica1dim}, one gets 
the following lemma (corresponding to Lemma~\ref{l.descont}):

\begin{lemm}
\label{l.WuQt}
Given any subinterval $I$ of
$D_t^+$ there are $x \in I$ and $i\geq 0$ with 
$\mfG_t(\mfR_t^i(x))=0$.
\end{lemm}

\subsection{Lateral Homoclinic classes. Proof of 
Theorem~\ref{t.latlat}}
\label{ss.classeshomoclinicasdaselano}

To prove Theorem~\ref{t.latlat} we proceed as
in Section~\ref{s.pontosnaoerrantes}.
After redefining vertical strips and segments and  using Lemma~\ref{l.WuQt},
one gets that,
for any
small $t>0$ and any
disk $\chi$ with 
$W^s(P,f_t)\pitchfork \chi\ne \emptyset$,
$W^{ss}(S,f_t)\pitchfork \chi\ne \emptyset$
(recall
Proposition~\ref{p.chi}).
The inclusion
$(H(P,f_t)\cup \Ups_t)\subset H^+(S,f_t)$
follows exactly as $(H(P,f_t)\cup \La_t)\subset H(Q,f_t)$
in the case of heterodimensional cycles.

\subsection{Proof of Theorem~\ref{t.HssubsetHp}:
the inclusion $H^+(S,f_{t})\subset H(P, f_t)$}
\label{ss.newingredients}

Consider the homoclinic point $x_t=(-1/2,t,0)$ of $P$ for $f_t$
and the fundamental domains
$\De^+_t(i) =
\Phi^{-i} (\De^+_t(0))$, $i\ge 0$,
where
$\De_t^+(0)=
[\Phi^{-1}(t), t]$.
Let
$\ka_t$
be the first $k \in \NN$ such that 
$\Phi^{k}(\De_t^+(0))\subset [1-t,1]$. 
Observing that, for small $t>0$, 
$|\Phi^{k_t}(\De_t^{+}(0))|\le t\, (1-\la)$ and
$|\De_t^+(0)|\ge \la \, t^2$, we get, 
using the Bounded Distortion Lemma~\ref{f.distorcao},
\begin{equation}
\label{f.kappa_t}
  (\Phi^{\ka_t})^\prime(x) <  \frac{1-\la}{\la \, t}\, e^{K},
\quad
\mbox{for all $x\in \De_t^+(0)$.}
\end{equation}
Denote by $\de_t^i$ the length of
$\De^+_t(i)$.
Since the derivative of $\Phi$ near $0$ is close to $1$
and strictly bigger than $1$ in $(0,t]$,
for small $t$, we have that
\begin{equation}
\label{f.9/10}
\displaystyle{\de_t^0 \geq \de_t^i \geq \frac{9\,\de_t^0}{10}},
\quad
i = 1,\dots,4.
\quad
\mbox{In particular,}
\quad
 \sum_{i = 0}^4 \de_t^i \in [4 \de_t^0, 5 \de_t^0].
\end{equation}

We now construct a family 
$\cH_t$
of homoclinic points of
$P$ for $f_t$ such that
the set $\{ y\colon \, (x,y,0)\in \cH_t\}$  
is dense in the fundamental domain $\De_t^+(0)$:

\begin{prop}
\label{p.sexsdensas}
For every small
$t>0$ there are sequences of homoclinic points of 
$P$ of the form 
$(b_{i_1,i_2,\dots, i_m,k},x_{i_1,i_2,\dots, i_m,k}, 0)_{k\in \NN^*}$,
$b_{i_1,i_2,\dots, i_m,k} \in [-1,0]$,
such that
\begin{description}
\item{{\bf (H1)}}
$x_{i_1,i_2,\dots, i_m,k}\in \bigcup_{i=0}^4 \De_t^+(i) = \De_t$,
\item{{\bf (H2)}}
$x_{i_1,i_2,\dots, i_m,k}\to 
x_{i_1,i_2,\dots, i_m}$ 
as $k\to \infty$,
\item{{\bf (H3)}}
$x_{i_1,i_2,\dots, i_m,0} <x_{i_1,i_2,\dots, (i_m-1)}$
for every $i_m\ge 1$,
\item{{\bf (H4)}} 
$\mbox{{\rm diam}}((x_{i_1,i_2,\dots, i_m,k})_k)\to 0$ as $m\to \infty$,
\item{{\bf (H5)}}
$(x_i)$ is increasing and $x_i\to t^-$ as $i\to \infty$,
\item{{\bf (H6)}}
$x_0\in \De_t^+(1)$ and  $x_0 \not\in \De_t^+(0)$.
\end{description}
\end{prop}

This proposition will be proved in Section~\ref{sss.sexsdensas}.
From the proposition one gets the following:

\begin{coro}
\label{c.chato}
The set
$\cH_t = \bigcup_{n,k \in \NN^*}(x_{i_1,i_2,\dots, i_k,n})$ 
contains a  dense subset of $\De^+_t(0)$.
\end{coro}

\begin{demo}
The proof of is identical to \cite[Lemma 4.1]{Dnon},
 but we repeat it here for
completeness.
Take any point $x\in \De_t^+(0)$. 
If
$x \in \cH_t$ there is nothing to prove.
Otherwise, 
by (H5) and (H6), there is 
$i_1 > 1$ with
$x_{(i_1-1)} < x < x_{i_1}.$ 
Analogously, 
by (H2) and (H3),  
there is $i_2>1$ with 
$x_{i_1,(i_2-1)} < x < x_{i_1,i_2}.$
Inductively,  using (H2) and (H3) as above,
we get a sequence $\{i_k\}$, $i_k > 1$,
such that, for all $k$, 
$
x_{i_1,\dots,(i_k-1)} < x < x_{i_1,\dots,i_k}.
$
Finally, from (H4), 
$\lim_{k \to \infty} x_{i_1,\dots,i_k} = x$, ending the proof of the
lemma.
\end{demo}

\subsubsection{Proof of Theorem~\ref{t.HssubsetHp}}
The deduction of Theorem~\ref{t.HssubsetHp} from
Corollary~\ref{c.chato} follows as in \cite[Section 5]{Dnon}.
For completeness, we sketch here this proof.
%
Consider any $w \in W^u(S,f_t)\pfrk W^{ss}(S,f_t)$.
By replacing $w$ by some iterate of it, we can assume that
$w=(x,0,0)$,  $|x|$ small.
We prove that, for  every $\ve>0$, 
the square
$S(\ve)=(x-\ve,x+\ve) \times (0,\ve)\times \{0\}$ 
contained in $W^s(P,f_t)$ transversely intersects 
$W^u(P,f_t)$.
This immediately implies
that $w\in H(P,f_t)$.
The configuration of the cycle and the $\la$-lemma
imply that there is $n(\ve)>0$ such that
$f_t^{-n(\ve)}(S(\ve))$
contains a disk 
$S^\prime(\ve)$
of the form
$$
S^\prime(\ve)= [-1,1]
\times (\bar y -\xi, \bar y +\xi)
\times \{\bar z\},
\quad
\mbox{$\bar y\in (1-t,1)$, $\bar z\in [-1,1]$ and small $\xi>0$.}
$$
Let $m\in \NN$ be such that
$\Phi^{-m}(\bar y)\in \De_t^+(0)$.
Thus $f_t^{-m}(S^\prime(\ve))$
contains the strip 
$$
\hat S(\ve)=
[-1,1] \times 
(\Phi^{-m}(\bar y-\xi), 
\Phi^{-m}(\bar y+\xi))
\times\{\la_u^{-m}\, \bar z\}
\subset
W^s(P,f_t).
$$
Since $\Phi^{-m}(\bar y)$ belongs to $\De_t^+(0)$,
Corollary~\ref{c.chato} implies that
$\hat S(\ve)$ meets $W^u(P,f_t)$. Thus $\hat S(\ve)$ contains 
a homoclinic point of $P$ and the same holds for $S(\ve)$.
\hfill
$\Box$

\subsubsection{Proof of 
Proposition~\ref{p.sexsdensas}: Sequences of homoclinic points:}
\label{sss.sexsdensas}

Consider the interval 
 $[1-\eta_t, 1]$, where
$\eta_t=\de_t^1+\de_t^0.$
Define $\alpha_t$ as the first natural number $\alpha$  with
$$
\Phi^{\ka_t+\alpha}(\De_t^+(0))\subset [1-\eta_t, 1].
$$
Observe that $|\De_t^+(0)| = \de_t^0<t^2$ 
and $\de_t^1<\de_t^0$. Thus, for small $t>0$, 
$\eta_t <  2\,\de_t^0< 2\, t^2<\la\, t$.
Since, by definition,
$\Phi^{\ka_t}(\Phi^{-1}(t))\in [1-t, 1-\la\, t]$ and
$(1-\la\,t)<(1-\eta_t)$,
we get that $\alpha_t\ge 1$ for every small
$t$.
Observe also that
$t^2 <\eta_t<2 \,t^2$,
where the first inequality follows from (\ref{f.9/10})
and
$\de_t^0>3t^2/4$ if $t$ is small enough.

\begin{lemm}
\label{l.lambdaalphat}
For every small $t>0$ it holds
$\la^{\alpha_t}\le  (2\, t)/\la.$
\end{lemm}

\begin{demo}
By definitions of $\ka_t$ and $\alpha_t$, $\Phi^{\ka_t}(\De_t^+(0))=
[1-e_t^-,1-e_t^+]$, where
$e_t^-\in [\la\, t, t]$, and
$\Phi^{\alpha_t}(1-e_t^-)\in [1-\eta_t, 1-\la\, \eta_t].$
Thus,
since
$\Phi$ is affine near $1$,
$\la^{\alpha_t}(e_t^-) \in (0,\eta_t].$
Thus, from $t^2<\eta_t<2\,t^2$ and 
$\la\, t \le e_t^-\le t$,
$$
\la^{\alpha_t} \leq \frac{\eta_t}{e_t^-} \leq 
\frac{2\, t}{\la},
$$
ending the proof of the lemma.
\end{demo}

Next contraction lemma is necessary for getting 
(H4) in  Proposition~\ref{p.sexsdensas} and 
along the inductive definition of the sequences
$(x_{i_1,i_2,\dots, i_m,k})_k$.

\begin{lemm}
\label{l.derivada<1/2}
$L = \max \{ (\Phi^{\ka_t+\alpha_t+j})^\prime (x);\,
x\in \cup_{i=0}^4 \De^+_t(i) \mbox{ and } j\ge 0\}< \frac{1}{2}.$
\end{lemm}

\begin{demo}
Since $\Phi$ is a contraction near $1$,
it is enough to compute the estimate when $j=0$.
We split the trajectory of a point $x\in  \De^+_t(i)$
going from $\De_t^+(i)$ to $[1-\eta_t,1)$ as follows:
{\bf (i)}
$i$ iterates,
$i\le 4$, for $x$  going  from $\De^+_t(i)$ to
$\De^+_t(0);$
{\bf (ii)}
$\ka_t$ iterates for $\Phi^i(x)$ 
going from $\De^+_t(0)$ to
$\Phi_t^{k_t}(\De^+_t(0))$; and 
{\bf (iii)}
$\alpha_t$ iterates  for  $\Phi^{\ka_t+i}(x)$ 
going from
$\Phi_t^{k_t}(\De^+_t(0))$
to
$ [1-\eta_t,1].$ 
This  construction involves
$(i + \ka_t + \alpha_t)$ iterations of $x$ by $\Phi$, that is,
we need to remove the last $i$  
iterations, 
corresponding to a contraction by $\la^i$.
We claim that
\begin{equation}
L
\leq 
\underbrace{((2\,t+1)^4)}_{(\mbox{{\bf a}})} 
\, 
\underbrace{
\frac{e^{K}\,(1 - \la)}{\la \, t} 
\,
\la^{\alpha_t}
}_{(\mbox{{\bf b}})}
 \,
\underbrace{\frac{1}{\la^4}}_{(\mbox{{\bf c}})},
\end{equation}
corresponding
(\mbox{{\bf a}}) to the expansion of the first
$i$ iterates by $\Phi$
(just observe that in $[0,t]$ the derivative of 
$\Phi$ is upper bounded by
$(2\,t+1)$  and that $i\le 4$),
(\mbox{{\bf b}})
to an upper bound of the derivative of
$\Phi^{\ka_t+\alpha_t}$, 
recall (\ref{f.kappa_t}), and
(\mbox{{\bf c}})  to the  $i$ ($i \leq 4$)
negative  iterates of $\Phi$ close to $1$.
By (SN), Lemma \ref{l.lambdaalphat}, and the fact that 
 ${(2\,t+1)}^4 < 2$ if $t>0$ is small, we get
$$
L
\leq
 ({(2\,t+1)}^4) \, 
\frac{e^{K}\,(1 - \la)}{\la\, t} \,
\frac{2\, t}{\la} \, \frac{1}{\la^4} 
= {(2\,t+1)}^4 \, \frac{2\,(1 - \la)\, e^{K}}{\la^6}
\leq 4\, e^{K}\,\frac{(1 - \la)}{\la^6} < \frac{1}{2},
$$
which ends the proof of the lemma.
\end{demo}

\paragraph{Construction of the sequences 
$(x_{i_1,i_2,\dots, i_m, k})$.} 
To construct  the 
sequences $(x_{i_1,i_2,\dots, i_m, k}),$
we need the following algorithm 
about the creation of homoclinic points,
which is a consequence of the definition of
the unfolding of the heterodimensional cycle.

\begin{algo}
\label{a.algo2}
Let $(x,y,0)$, $y\in [0,t]$, be a homoclinic point of
$P$ (for $f_t$) such that
$\{(x,y)\}\times [-1,1]\subset W^u(P,f_t)$.
Then, for every $m$ with
$\Phi^m(y)\in (1-t,1)$,
 there is a homoclinic point
of $P$ of the form
$(\bar x, \Phi^m(y)+t-1,0)$
such that
$\{(\bar x, \Phi^m(y)+t-1)\}\times [-1,1]\subset
W^u(P,f_t)$.
\end{algo}

Consider the homoclinic point
$(-1/2,t,0)$ of $P$ for $f_t$
(satisfying Algorithm~\ref{a.algo2}) and
the sequences
$(y_i)_{i \in \NN^*}$ and 
$(x_i)_{i \in \NN^*}$ defined by
$$
y_i= \Phi^{\ka_t+ \alpha_t + i}(t) \quad 
\mbox{and}
\quad
x_i=(t-1)+y_i,
\quad
y_i\to 1
\quad
\mbox{and}
\quad
x_i \to t.
$$
Observe that, for each $i\ge 0$, there is a homoclinic
point $(b_i,x_i,0)$ of $P$ verifying Algorithm~\ref{a.algo2}.
Also,
by the definitions of
$\alpha_t$ and $\ka_t$, 
$y_i\in [1-\eta_t,1]$ for all $i\ge 0$.
Thus, since $\eta_t=\de_t^0+\de_t^1$, one has
\begin{equation}
\label{f.x_iinDt}
x_i\in [t- \eta_t,t] = [t-(\de_t^1+\de_t^0),t]
 = \De_t^+(1)\cup \De_t^+(0).
\end{equation}

\begin{lemm}
\label{l.h6} 
The sequence
$(x_i)_{i \in \NN^*}$ verifies (H5) and (H6).
\end{lemm}

\begin{demo}
Condition (H5) follows by definition. 
To get (H6), i.e. $x_0 =  y_0 +(1-t) \in (\De_t^+(1)\setminus \De_t^+(0))$,
observe that,
by  construction, 
$x_0\in [t-\eta_t, t- \la \, \eta_t]$
and
$\De_t^+(1)= [t-\eta_t, t-\de_t^0]$. Thus, 
by (\ref{f.x_iinDt}), it is enough 
to check that
$\la\, \eta_t
=\la\,(\de_t^0+\de_t^1)>\de_t^0$.
This inequality folllows from $\de_t^0 >\de_t^1$, (\ref{f.9/10}) and (SN), 
observing that
$$
\frac{\de_t^0}{\de_t^0+\de_t^1}< \frac{\de_t^0}{2\, \de_t^1}<
\frac{\de_t^0}{2\, (9/10) \, \de_t^0}=\frac{10}{18}<
\frac{2}{3}
<\la.
$$
The proof of the lemma is now complete.
\end{demo}

We now proceed with the construction  of
the sequences in Proposition~\ref{p.sexsdensas}.
For each
$j\in \NN^*$, 
define the sequences
$(y_{j,i})_{i \in \NN^*}$ and 
$(x_{j,i})_{i \in \NN^*}$ as follows,
$$
y_{j,i}=\Phi^{\ka_t+\alpha_t+j}(x_i)
\quad
\mbox{and}
\quad
x_{j,i}=(t-1) + y_{j,i}.
$$
We claim that
$y_{j,i} \to y_j$ and, consequently,
$x_{j,i} \to x_j$, as $i \to \infty$. For that just observe that
$\lim_{i\to\infty}x_i=t$, thus, by
continuity,
$
\lim_{i\to\infty}y_{j,i} = 
\lim_{i\to\infty}\Phi^{\ka_t+\alpha_t+j}(x_i)=
\Phi^{\ka_t+\alpha_t+j}(t)=y_{j}.$

\begin{lemm}
The points $(x_{j,i})_i$ belong to
$\cup_{i=0}^4 \De_t^+(i)$
for all  $i, j\in\NN \cup \{0\}$. 
\label{l.step2}
\end{lemm}

\begin{demo}
By construction, the sequences $(x_{j,i})_i$ are increasing, thus
it is enough to prove that $x_{0,0}\in  \cup_{i=0}^4 \De_t^+(i)$.
Consider the diameter  $d_0=(t-x_0)$ of $(x_i)_{i\in \NN^*}$.
By (\ref{f.x_iinDt}),
 $d_0<\de_t^0+\de_t^1$.
Let 
$d_1=|x_0-x_{0,0}|$
be the diameter of 
the sequence $(x_{0,i}=\Phi^{\ka_t+\alpha_t}(x_i)+(t-1))_i$, 
which is equal to the diameter of $(\Phi^{\ka_t+\alpha_t}(x_i))_i$.
Therefore,
since $(x_i)_i\subset \De_t^+(0)\cup \De_t^+(1)$,
by Lemma~\ref{l.derivada<1/2},
the diameter
$d_1$  is bounded by
\begin{equation}  
\label{e.step2}
d_1\le L \, d_0< d_0/2< (\delta^0_t+\delta^1_t)/2< 
(2\, \delta^0_t)/2<
\delta^0_t.
\end{equation}
Since $x_0\in \De_t^+(1)$ (Lemma~\ref{l.h6}),
to prove that  $x_{0,0}\in  \cup_{i=0}^4 \De_t^+(i),$
it is enough to see that
$$
x_{0,0}
=x_0-d_1>
x_0-(\de_t^2+\de_t^3+\de_t^4)
\Longleftrightarrow
d_1<(\de_t^2+\de_t^3+\de_t^4),
$$
which immediately follows from
$\de_t^2+\de_t^3+\de_t^4> \de_t^0>d_1$, the first inequality being consequence
of
$\de_t^i>(9\, \de_t^0)/(10)$, see (\ref{f.9/10}), and the last from 
(\ref{e.step2}).
This ends the proof of the lemma.  
\end{demo}

Suppose now inductively defined 
sequences
 $(y_{i_i,i_2,\dots,i_m,i})_{i\in \NN^*}$ and 
$(x_{i_i,i_2,\dots,i_m,i})_{i\in\NN^*}$
by
$$
y_{i_1,i_2,\dots,i_m,i}=\Phi^{\ka_t+\alpha_t+i_1}
(x_{i_2,\dots,i_m,i})
\quad
\mbox{and}
\quad
x_{i_1,i_2,\dots,i_m,i}=
(t-1) + y_{i_1,i_2,\dots,i_m,i},
$$
satisfying conditions
(H1), (H2), (H3) and
\begin{description}
\item
{{\bf (H2b)}}
$(y_{i_1,i_2,\dots,i_m,i})_i \to y_{i_i,i_2,\dots,i_m}$
as $i\to \infty$,
\item
{{\bf (H3b)}}
$y_{i_1,i_2,\dots,i_m,0}<  y_{i_i,i_2,\dots,(i_m-1)}$
for all  $i_m\ge 1$.
\item{{\bf (H4b)}}
Let  $d_m$, $m\ge 0$, be the  diameter of the  sequence 
$(x_{\underbrace{0,\dots,0,i}_{m \, 0's}})_{i}$. 
Then  $d_m\le  (d_{m-1})/2$.
\end{description}

Observe that
(H2) and (H2b) (resp. (H3) and (H3b)) are equivalent.
Notice that, for $m=1$,
(H1) follows from Lemma~\ref{l.step2}, (H2) and (H2b) from definition, 
and (H4b) from (\ref{e.step2}).
To check (H3b), $y_{i,0}<y_{i-1}$, for every $i\ge 1$, 
recall that,
by Lemma~\ref{l.h6},
$x_0<\Phi^{-1}(t)<t$, thus
$$
y_{i-1}=\Phi^{\ka_t+\alpha_t+i-1}(t)=
\Phi^{\ka_t+\alpha_t+i}(\Phi^{-1}(t))>
\Phi^{\ka_t+\alpha_t+i}(x_0)=y_{i,0}.
$$

For  simplicity we say that
the sequences 
$(z_{i_i,i_2,\dots,i_m,i})_{i\in\NN^*}$, $z=x,y$, are
of generation $m$.

\begin{lemm}
\label{l.i4h4}
Property (H4b) implies (H4) 
in  Proposition~\ref{p.sexsdensas}.
\end{lemm}

\begin{demo}
By construction,
the diameters of the  sequence of generation $m$ are
bounded by the diameter $d_m$ of
$(x_{0,\dots,0,i})_{i\in\NN^*}$.  Thus, inductively, 
$d_m\le (1/2) \, d_{m-1}\le (1/2)^{m} \, d_0,$ so $d_m\to 0$.
\end{demo}

Keeping in mind Lemmas~\ref{l.h6} and \ref{l.i4h4}, in order 
to prove Proposition~\ref{p.sexsdensas} it suffices to
see that the sequences above verify (H1), (H2b), (H3b) and (H4b). 
We argue inductively
on the generation of the sequences and assume satisfied these conditions
for sequences of generation less than or equal to
$m$.
To verify 
(H2b) 
for the sequences of generation $m+1$
note that, by induction,
$(y_{i_1,i_2,\dots,i_m,i})_i \to y_{i_1,i_2,\dots,i_m}$. Thus,
by continuity of $\Phi$ and by definition,
$$
(y_{j,i_1,i_2,\dots,i_m,i})_i=
(\Phi^{\ka_t+\alpha_t+j}(x_{i_1,i_2,\dots,i_m,i}))_i 
\to 
\Phi^{\ka_t+\alpha_t+j}(x_{i_1,i_2,\dots,i_m})
=y_{j,i_1,i_2,\dots,i_m}.
$$
To prove   (H3b) observe that, by induction,
$
x_{i_1,i_2,\dots,i_m,0}
<x_{i_1,i_2,\dots,i_m-1}.
$
Thus,  since $\Phi$ is increasing,
$$
y_{j,i_1,i_2,\dots,i_m,0}=
\Phi^{\ka_t+\alpha_t+j}(x_{i_1,i_2,\dots,i_m,0})
<
\Phi^{\ka_t+\alpha_t+j}
(x_{i_1,i_2,\dots,i_m-1})=
y_{j,i_1,i_2,\dots,i_m-1}.
$$
To check (H4b) 
observe that, by the induction hypotheses (H1), 
$x_{0,\dots,0,i}\in  \cup_{i=0}^4 \De_t^+(i)$, 
Lemma~\ref{l.derivada<1/2} 
and the fact that the sequences 
$(y_{0,\dots,0,i})$ and  
$(x_{0,\dots,0,i})$ have the same diameter
imply that
$$
d_{m+1} =\mbox{diam}(\underbrace{
                                 (y_{0,\dots,0,i})_i}_{(m+1)\, 0's})
=
\mbox{diam}((\Phi^{\ka_t+\alpha_t}
                    (\underbrace{(x_{0,\dots,0,i})_i)}_{m \, 0's})
=L \,
\mbox{diam}
((x_{0,\dots,0,i}))_i) = L \, d_m \le \frac{d_m}{2},
$$
which ends the proof of (H4b).

Finally, to get (H1) for sequences of generation $(m+1)$ 
it is enough to see that, for every $m$,
$
\sum_{i=0}^m d_m < 4 \de_t^0 < \sum_{i=0}^{4} \de_t^i,
$
where the last inequality follows from (\ref{f.9/10}).
By induction and
Lemma~\ref{l.h6},
which implies that $d_0\le \de_t^0+\de_t^1$,
we have
$$
\sum_{i=0}^m d_m \le
\sum_{i=0}^{m-1} (1/2)^i \, d_0 \le 
\sum_{i=0}^{m-1} (1/2)^i \, 
 (\de_t^0+\de_t^1) \le 
\frac{1}{1-2}\, (\de_t^0+\de_t^1)<
{2}(\de_t^0+\de_t^1)< 4\de_t^0,
$$
ending the proof of our claim

The construction of the sequences  
$(x_{i_i,i_2,\dots,i_m,i})_{i\in\NN^*}$ of 
Proposition~\ref{p.sexsdensas} is now complete.

\section{Homoclinic classes before collapsing
the saddles $S_s^+$ and $S^-_s$}
\label{s.lateralhomoclinic}

We now return to the two parameter family $f_{t,s}$
in Section~\ref{s.selano}.
Observe that
$S_s^+=(0,\sqrt{s},0)$, $s>0$, is a fixed point of index two
of $f_{t,s}$  (any $t>0$)
and that
$f_{\sqrt{s},s}$ has a heterodimensional cycle
associated to $S^+_s$ and $P$:
\begin{itemize}
\item
$W^u(S_s^+,f_{\sqrt{s},s})$
 meets transversely
$W^s(P,f_{\sqrt{s},s})$ 
throughout the segment $\{0\}\times (\sqrt{s},1)\times \{0\}$,
\item
$W^u(P,f_{\sqrt{s},s})$ meets quasi-transversely  $W^s(S^+_s,f_{\sqrt{s},s})$
along the orbit of
$(-1/2,\sqrt{s},0)$
(just observe that
$[-1,1]\times \{(\sqrt{s},0)\}
\subset W^s(S_s^+,f_{\sqrt{s},s})$
and that $(-1/2,\sqrt{s},0)\in W^u(P,f_{\sqrt{s},s})$).
\end{itemize}

In what follows we assume that the saddle-node arc $\Phi_s$ 
verifies condition (SN) in Section~\ref{s.selano}.

\begin{theor}
\label{t.lats}
There exist a  small $s_0>0$ and a strictly positive map 
$\tau$
defined on $(0,s_0)$  such that,
for every $s\in (0,s_0)$ and $t\in (\sqrt{s},\sqrt{s}+\tau(s))$,
\begin{itemize}
\item
$H(P,f_{t,s})=
H(S^+_{s},f_{t,s})$, and 
\item
there is a neighborhood $W_s$ of the
cycle of $f_{\sqrt{s},s}$ 
(associated to $P$ and $S^+_s$)
such that the maximal invariant
set $\La_{t,s}$ of $f_{t,s}$ in $W_s$ is equal to $H(S^+_s,f_{t,s})$.
\end{itemize}
\end{theor}

The proofs of the  inclusion    
$H(P,f_{t,s})\subset
H(S_{s},f_{t,s})$ and the second part of
the  theorem
follow as in  Theorem~\ref{t.lat}.
So we just sketch these proofs.
For a fixed $s>0$ and $t>\sqrt{s}$, $t$ close to 
$\sqrt{s}$,
$t=\sqrt{s}+\tau$,
consider the scaled fundamental domains $D_{t,s}^\pm$
of $\Phi_s$,
$$
D_{t,s}^-=[1-(t-\sqrt{s}), 1-\la\, (t-\sqrt{s})]=
[1-\tau, 1-\la \, \tau]
$$ 
and 
$D_{t,s}^+$ defined as the first backward iterate of
$D_{t,s}^-$ in 
$[\sqrt{s}, \sqrt{s}+\tau]$. Let 
$\Phi_s^{k_{t,s}}(D_{t,s}^+)= 
D_{t,s}^-$, where $k_{t,s}\in \NN$.
These domains play the role of  
$D_t^\pm$ in Section~~\ref{s.dinamica1dim}.
Observe that 
$$
\ell (t,s)=\frac{|D_{t,s}^-|}{|D_{t,s}^+|}
\ge
\frac{\tau\,(1-\la)}
{(\sqrt{s} +\tau)^2-s}
=
\frac{\tau\,(1-\la)}
{\tau\, (2\sqrt{s}+\tau)}
\ge
\frac{\tau\,(1-\la)}{\tau\, (2\sqrt{s}+\tau)}=
\frac{1-\la}{2\, \sqrt{s}+\tau}.
$$
By shrinking $s$, we can assume that
$|\Phi^{\prime\prime}_s(x)|/|\Phi^\prime_s(x)|<K$ for all $x\in [-1,2]$
($K$ as in Lemma~\ref{f.distorcao}).
Thus, there is  $s_0>0$ and a map $\tau\colon (0,s_0)\to
\mathbb{R}^+$
such that 
\begin{equation}
\label{e.expansion}
\ell(t,s) \, e^{-K}>2,
\quad
\mbox{for all $s\in (0,s_0)$ and $t\in (\sqrt{s},\sqrt{s}+\tau (s))$}.
\end{equation}
Exactly as in Section~\ref{s.dinamica1dim},
for  $s\in (0,s_0)$ and $t\in (\sqrt{s},\sqrt{s}+\tau (s))$,
we define 
maps
$$
\begin{array}{ll}
&
T_{t,s}\colon  D_{t,s}^+ \to  D_{t,s}^-,
\quad T_{t,s}(x)=\Phi_s^{k_{t,s}}(x),\\
&
G_{t,s}\colon  D_{t,s}^- \to  [\sqrt{s},
\sqrt{s}+\tau \,  (1-\la)],
\quad
G_{t,s}(x)=T_{t,s}(x)+(t-1),
\\
&
R_{t,s}\colon D_{t,s}^+\to D_{t,s}^+,
\quad
R_{t,s}(x)=
\Phi_{s}^{i(x)} (G_{t,s}(x)),
\end{array}
$$
where, as in Section~\ref{s.dinamica1dim},
$i(x)$ is the first forward iterate of 
$G_{t,s}(x)$ by $\Phi_s$ in $D_{t,s}^+$.

As in the 
proof of Lemma~\ref{l.mfTtexpansora},
the Bounded Distortion Lemma~\ref{f.distorcao}
and equation (\ref{e.expansion})
imply that the maps  $T_{t,s}$, $G_{t,s}$ and $R_{t,s}$
are all uniformly expanding.
The inclusions 
$H(P,f_{t,s})\subset
H(S^+_{s},f_{t,s})$ and $\La_{t,s}= H(S_s^+,f_{t,s})$
now follow as in the proof of Theorem~\ref{t.lat}.

To prove $H(P,f_{t,s})\subset
H(S_{s},f_{t,s})$, recall 
that Theorem~\ref{t.HssubsetHp} 
gives small $\bar t>0$ with
$H^+(S,f_{t,0})\subset H(P,f_{t,0})$ for all $t\in (0,\bar t)$.
The proof of this inclusion follows after 
constructing the  sequences of
homoclinic points of $P$ in
Proposition~\ref{p.sexsdensas}.
The proof of such a proposition only involves distortion 
control of the saddle-node map in  
$[0,1]$ 
and the definition of contracting itineraries,
(Lemma~\ref{l.derivada<1/2}). 
Clearly, these properties also hold after replacing, for small
positive $s$, 
the saddle-node
$S$ by the hyperbolic point $S^+_s$ and considering the
restriction of the saddle-node map to 
$[\sqrt{s},1]$. The sketch of the proof of Theorem~\ref{t.lats} is
now complete.

\section{Collision, explosion and collapse of homoclinic classes}
\label{s.colapso}

In this section we prove Theorems~\ref{t.main} and \ref{t.saturated}. 
Consider a two
parameter family of diffeomorphism $f_{t,s}$  
locally defined as follows:

\medskip
\noindent
{\bf Partially hyperbolic local dynamics:}
\begin{itemize}
\item
In the set
$C=[-1,1]\times [-2,2]\times [-1,1]$,  
$f_{t,s}(x,y,z)=
(\la^s x, \Psi_s(y),
\la^u z)$,
where
$0<\la^s<1<\la^u$ and
$\Psi_s\colon [-2,2]\to (-3,2)$ 
is a strictly increasing
map
such that
$\la^s<d_m<\Psi_s^\prime(x)<d_M<\la^u$
for all $x\in [0,1]$ and small $|s|$.
\item
$\Psi_s(1)=1$ and $\Psi_s(-1)=-1$ for all $s$. Moreover,
there exists $\delta>0$ such that,
in the intervals $[-1-\delta,-1+\delta]$ and 
$[1-\delta,1+\delta]$,
$\Psi_s$   
is affine and independent of $s$.
Furthermore, $\Psi_s^\prime(-1)=\beta$ 
$\Psi_s^\prime(1)=\la$,
where $0<\lambda<1<\beta$.
\item
There is $\delta>0$ such that the restriction
of $\Psi_s$ to $[-\delta,\delta]$ is
of the form
$\Psi_s(x)=x+x^2-s$.
\item
For $s<0$, $\Psi_s$ has (exactly) two fixed points
in $[-2,2]$ (the points $\pm 1$), for $s=0$,
$\Psi_0$ has (exactly) three fixed points
($\pm 1$ and $0$), and, for $s>0$, $\Phi_s$ has 
(exactly) four fixed points ($\pm 1$ and $\pm \sqrt{s}$).
\end{itemize}

We let $P=(0,1,0)$, $Q=(0,-1,0)$, $S=(0,0,0)$ and,
for positive $s$, $S^\pm_s=(0,\pm \sqrt{s},0)$, the
fixed points of $f_{t,s}$ in $C$.

\medskip
\noindent
{\bf Existence and unfolding of cycles:}
\begin{itemize}
\item
There are $k_0\in \NN$ 
 and small neighborhoods of 
$Q$, $S$ and $P$
such that, for each small $|s|$, 
in these neighborhoods  $f_{t,s}^{k_0}$
is the translation
$
f^{k_0}_{t,s}(x,y,z)=
(x-1/2, y-1+t, z+1/2),
$
recall the definition of $f_t$ in Section~\ref{s.heterodimensional}.
\end{itemize}
As in previous sections, we have:
\begin{itemize}
\item
for $(t,s)=(0,0)$, $f_{0,0}$ exhibits a pair of  saddle-node
heterodimensional cycles, 
associated to 
$P$ and $S$ 
($W^u(P, f_{0,0})$  meets transversely $W^{ss}(S, f_{0,0})$)
and to
$Q$ and $S$ (here $W^s(Q, f_{0,0})$  intersects $W^{uu}(S, f_{0,0})$).
Just observe that $[-1,1]\times \{(0,0)\}\subset
W^{ss}(S,f_{0,0})$ and   $\{(0,0)\}\times [-1,1]\subset
W^{uu}(S,f_{0,0})$.
\item
For small $|t|$ and $s<0$, 
the homoclinic classes of $P$ and $Q$ are both
nontrivial: 
notice that, for negative $s$,  
$[-1,1]\times (-1,2)\times \{0\}\subset
W^s(P,f_{t,s})$ and 
$\{0\}\times (-1,2)\times [-1,1]\subset
W^u(Q,f_{t,s})$, thus  
$(-1/2,t,0)$ and
$(-1/2,-1,0)$ are  homoclinic points of $P$ and $Q$,
respectively.  
\item
$f_{\sqrt{s},s}$
has a pair of heterodimensional cycles associated to
$P$ and $S^+_s$ and to $Q$ and $S^-_s$
(this is obtained exactly as in Section~\ref{s.lateralhomoclinic}).
\end{itemize}

As in Section~\ref{s.selano}, we
assume the following distortion property
(similar to (SN)) involving the saddle-node $S$ and
the saddles $P$ and $Q$:
Let $K$ be an upper bound for
$|\Psi_s''(x)|/|\Psi_s'(x)|$, for small $|s|$ and $x\in 
[-1,1]$,
\begin{description} 
\item
{\bf (DS)}
$
\max
\left\{
\displaystyle{   
\frac{4 \, e^{K}\, (1-\la)}{\la^6},
\frac{4 \, e^{K}\, (1-\be^{-1})}{\be^{-6}}}
\right\}
<
\displaystyle{ \frac{1}{2}}$,
where
$2/3<\la < 1< \beta<3/2$.
\end{description}

\subsection{Dynamics before collapsing  the saddles $S_s^+$ and $S_s^-$}
\label{ss.beforecollapse}

Theorems~\ref{t.latlat} and \ref{t.HssubsetHp} 
give the existence of a small positive  $\bar t$ such that
$H(P,f_{\bar t,0})=H^+(S,f_{\bar t,0})$ and
$H(Q,f_{\bar t,0})=
H^-(S,f_{\bar t,0})$.
The proof of these identities only involve
the following ingredients:
\begin{itemize}
\item The ratio between the lengths of
the scaled fundamental domains at the hyperbolic point and
at the saddle-node, and that such a ratio is arbitrarily large
(for the inclusions 
$H(P,f_{t_0,0})\subset H^+(S,f_{t,0})$ and
$H(Q,f_{t,0})\subset H^-(S,f_{t,0})$).
\item
The inclusion $H^+(S,f_{t,0})\subset  
H(P,f_{t_0,0})$ (resp. $H^-(S,f_{t,0})\subset  
H(Q,f_{t_0,0})$) is obtained by constructing sequences of
homoclinic points of $P$ (resp. $Q$) verifying 
Proposition~\ref{p.sexsdensas}.
The proof of such a proposition only involves 
distortion control of the saddle-node map in  
$[0,1]$ (resp. $[-1,0]$)
and construction of contracting itineraries
(Lemma~\ref{l.derivada<1/2}). 
\end{itemize}
Clearly, these properties also hold after replacing, for small
positive $s$, 
the saddle-node
$S$ by the hyperbolic point $S^+_s$ (considering the
restriction of $\Psi_s$ to 
$[\sqrt{s},1]$) and 
and the saddle-node by the point $S^-_s$
(considering the
restriction of $\Psi_s$ to
$[-1, -\sqrt{s}]$).
In this way, we get:

\begin{theor}
\label{t.iguais}
If $t>0$ is small, there is a small $s(t)>0$
such that
$H(P,f_{t,s})=
H(S_s^+,f_{t,s})$ and
$H(Q,f_{t,s})=
H(S_s^-,f_{t,s})$ for all $s\in (0,s(t))$.
\end{theor}

\subsection{Dynamics after collapsing the saddles $S_s^+$ and $S_s^-$}
\label{ss.aftercollapse}

\begin{theor}
\label{t.HP=HQ}
For every small  $t>0$ and $s<0$,
$H(P,f_{t,s})= H(Q,f_{t,s}).$
\end{theor} 

This theorem follows by using some of the ideas  in 
Theorem~\ref{t.lat}. 
We will prove the inclusion 
$H(P,f_{t,s})\subset H(Q,f_{t,s})$,  
($H(Q,f_{t,s})\subset H(P,f_{t,s})$  
similarly follows by considering
$f_{t,s}^{-1}$).
As in Sections~\ref{s.dinamica1dim}
and \ref{ss.1dim},
we define transitions $\cT_{t,s}$ and return maps
$\cR_{t,s}$ as follows.
Consider the fundamental domain 
$D_t^- = [1-3\, t, 1-3\, \la \, t]$
of $\Psi_s$ and define $k_{t,s}$ as the first $k\in \NN$ such that 
$\Psi_s(-t) \in \Psi_s^{-k}(D_t^-)$.
We let
$D_{t,s}^+ = \Psi_s^{-k_{t,s}}(D_t^-)$
and define the maps
$$
\cT_{t,s}\colon  D_{t,s}^+ \to  D_{t}^-,
\,\,\,     T_{t,s}(x)=\Psi_s^{k_{t,s}}(x),\qquad
\cG_{t,s}\colon  D_{t,s}^+ \to  [-2t, -t],
\,\,\,
\cG_{t,s}(x)=\cT_{t,s}(x)+t-1.
$$
Observe that,
by definition,
$\cG_{t,s}(D_{t,s}^+)\subset
[-2t,t(1-3\,\la)]\subset [-2t,-t]$
(recall $\la>2/3$).

Observing that, by definition of $\Psi_s$, $|D_{t,s}^+|\leq t^2 +s$,
the Bounded Distortion Lemma~\ref{f.distorcao}
and
$|D_t^-|=3\, \la\, (1-t)$
immediately give that
$$
\cT_{t,s}^\prime(x)
\ge
(e^{-K \, 2})\,
\frac{|D_{t}^-|}{|D_{t,s}^+|}
\ge  
(e^{-K\, 2})\, \frac{3\, \la(1-t)}{t^2+s}.
$$
 This inequality immediately implies the following:

\begin{lemm}
\label{l.28}
The maps 
$\cT_{t,s}$ and $\cG_{t,s}$ are $61$-expanding for all
small $t>0$ and $s<0$.
\end{lemm}

Since, by definition of $D_{t,s}$, $(-t)$ is at the left of
$D_{t,s}^+$, then, for each $x\in D_{t,s}^+$, there is a first $i(x)\ge 0$
such that
$\Psi_s^{i(x)}(\cG_{t,s}(x))\in D_{t,s}^+$.
We now define the return map $\cR_{t,s}$ by
$$
\cR_{t,s}
\colon D_{t,s}^+\to D_{t,s}^+,
\quad
\cR_{t,s}(x)=
\Psi_{s}^{i(x)} \cG_{t,s}(x).
$$

\begin{lemm}
\label{l.cRtexpansora}
The map  $\cR_{t,s}$
is $3$-expanding for all small $t>0$ and $s<0$.
\end{lemm}

Observe that, contrary what happens with the return maps
$R_t$ and $\mfR_t$ in Sections~\ref{s.dinamica1dim}
and \ref{ss.1dim}, 
the expansion for $\cR_{t,s}$ does not follow immediately 
from the expansion of
 $\cT_{t,s}$: the 
$i(x)$ iterates  by $\Psi_s$ at the left of $0$
contribute with a small contraction. 

\begin{demo}
By Lemma~\ref{l.28}, 
it is enough to verify that the contraction 
introduced by $\Psi_s^{i(x)}$ is at most
$1/20$.
For that recall that $\cG_{t,s}(D_{t,s}^+) \subset 
[-2\,t, -t]$ and observe that
 $i(x)\le i_0$, where $i_0$ is the first natural number with 
$\Psi^{i_0}(-2t)\in D_{t,s}^+$.
Write
$D_{t,s}^+(j)=
\Psi_s^{-j}(D_{t,s}^+)$
and note that, for every $x\in D_{t,s}^+$, 
$|D_{t,s}^+(i(x))|\le 
|D_{t,s}^+(i_0)|<
9\,t^2 +s$.
For each $i\in \{0,\dots, i_0\}$,
there is
$z_i \in D_{t,s}^+(i)$
such that
$$
({\Psi_s^{i}})^\prime(z_i)=
\frac{|D_{t,s}^+|}{|D_{t,s}^+(i)|}
\geq
\frac{|D_{t,s}^+|}{|D_{t,s}^+(i_0)|}
\geq
\frac{|D_{t,s}^+|}{9\,t^2+s}.
$$
Using the arguments
in the Bounded Distortion Lemma~\ref{f.distorcao}, 
we have that, for every $x \in D_{t,s}^+(i)$,
$$
({\Psi_s^{i}})^\prime(x)
\geq
\frac{e^{-\,K\, L_t}\,|D_{t,s}^+|}{9\,t^2+s},
$$
where $K>0$ is an
upper bound for  
$|\Psi_s''(x)|/|\Psi_s'(x)|$, $x\in [-1,1]$ and  $|s|$ small, and
$L_t=4\,t$ is the length of
$[-3\,t, t]$. 
Finally, observing that 
$|D_{t,s}^+| >(t^2)/2$,
$$
({\Psi_s^{i(x)}})^\prime(x)
\geq
\frac{e^{-K\, 4\,t}\,|D_{t,s}^+|}{9\,t^2+s} > 
\frac{e^{-K\, 4\, t }\,t^2}{18\,t^2+2s} >\frac{1}{20},
$$
for every small $t$ and $s$,
ending the proof of the lemma.
\end{demo}

Let  $D_{t,s}^+=[e_{t,s}^-,e_{t,s}^+]$. 
Observe that
there exist $i_1$ and $i_2\in \NN$ such that
$i(x)\in [i_1, i_2]$ for all $x\in D_{t,s}^+$ and, for each 
$i \in [i_1, i_2-1]$,
there is $d_i\in D_{t,s}^+$ with
$\Psi_s (\cG_{t,s}(d_i))=e_{t,s}^-.$ 
As in Section~\ref{s.dinamica1dim} and by definition, 
$d_{i_2-1}< d_{i_2-2}< \cdots < d_{i_1}$ and 
the points $d_i$ are the discontinuities of
$\cR_{t,s}$. Moreover, for each $i_1\le i<i_2$,
$\cG_{t,s}((d_{i+1}, d_i))=\mbox{int}(D_{t,s}^+)$
and $\cG_{t,s}$ is increasing in $(d_{i+1}, d_i)$.
Write $[d_{i+1}, d_i]= I_i$, $i_1 \ge i \ge i_2-2$,
$I_{i_2-1} = [e^-_{t,s}, d_{i_2-1}]$,  and
$I_{i_1} = [d_{i_1},e^+_{t,s}]$.
We continuously extend
$\cR_{t,s}$,
to the closed intervals $I_i$
(so $\cR_{t,s}$ is bivaluated at any $d_i$).

\begin{lemm}
\label{l.J}
Given any subinterval $J$ of $D_{t,s}^+$, there is 
$m\ge 0$ such that
$\cR_{t,s}^m(J)=D_{t,s}^+$.
\end{lemm}

\begin{demo}
It is enough to see
$\cR_{t,s}^m(J)$ contains an interval $[d_{i+1},d_i]$, 
$i_1\le i \le  i_2-2$, for some $m\in \NN^*$.
Write $J = J_0$.
If ${\cR_{t,s}}(J_0)$ contains two discontinuities we are done.
Otherwise, 
 ${\cR_{t,s}}(J_0)\subset I_i$  for some $i$ and  we let  
$J_1=\cR_{t,s}(J_0)$. By Lemma~\ref{l.cRtexpansora}, $|J_1|>3\, |J_0|$.
Finally,
if
${\cR_{t,s}}(J_0)$ just contains one discontinuity, say
$d_i$, then ${\cR_{t,s}}(J_0)$ only meets $I_i$ and $I_{i-1}$.
Write
$J_1^-={\cR_{t,s}}(J_0)\cap I_{i-1}$ and
 $J_1^+={\cR_{t,s}}(J_0)\cap I_i$, and let
$J_1$ be the biggest $J_1^\pm$.  
By Lemma~\ref{l.cRtexpansora}, $|J_1|>(3/2)\, |J_0|$. 
We now inductively define
intervals $J_i$ contained in the forward
orbit of $J_0$ such that 
either $J_{i+1}$ contains two discontinuities or
$|J_{i+1}|\ge (3/2) |J_i|$.  
Since the size of $D_{t,s}^+$ is finite,
at some step we get a first interval $J_m$ containing two discontinuities.
\end{demo}

Lemma~\ref{l.J} is the main step to prove 
Theorem~\ref{t.HP=HQ}, whose proof follows arguing as
in the proof of Theorem~\ref{t.lat} after the following considerations:

\begin{prop}
\label{p.faixavert2}
Let $\chi\subset
C=[-1,1]\times [-1,2]\times [-1,1]$, $x\in [-1,1]$
be a set of the form
$\chi=\{x\}\times A$, where $A$ is a disk of $\mathbb{R}^2$
whose interior contains 
a point of the form
$(y,0)$, with
$y\in (-1,2)$.
Then $\chi$ 
 intersects transversely $W^s(Q,f_{t,s})$.
\end{prop}

This result exactly corresponds
to Proposition~\ref{p.chip}. After proving it, 
Theorem~\ref{t.HP=HQ} is deduced as  follows.
As in  Section~\ref{s.pontosnaoerrantes}, 
Proposition~\ref{p.faixavert2}
implies that 
$W^s(Q,f_{t,s})$ meets transversely 
every
two-disk $\chi$ with 
$W^s(P,f_{t,s})\pitchfork \chi\ne \emptyset$.
Finally, arguing exactly as in 
Section~\ref{ss.homoclinicclasses},
we get      
$H(P,f_{t,s})\subset
H(Q,f_{t,s})$.

We now go to the details of the proof of Proposition~\ref{p.faixavert2}.
The first step is the following (corresponding to
Lemma~\ref{l.Wuptsegtovert}
in the proof of Propositions~\ref{p.chi} and \ref{p.chip}).

\begin{lemm}
\label{l.segto}
The manifold 
$W^u(P,f_{t,s})$ contains a vertical segment of the form
$\{(x,y)\}\times [-1,1]$, $x\in [-1,1]$ and $y\in D_{t,s}^+$. 
\end{lemm}

\begin{demo}
The proof follows as in
Lemma~\ref{l.Wuptsegtovert}, so we just sketch it.
Take the homoclinic point  $x_t=(-1/2,t,0)$ of $P$ for
$f_{t,s}$ and 
let  $r\in \NN$ be such that  
$\Psi^r_s(t)\in D_t^-$.
By construction,  
$$
\{(-1/2+ \la_s^r\, (-1/2),\Psi^r_s(t)+t-1)\}\times [-1,1]
\subset
W^u(P,f_{t,s}).
$$
Let $i=i(\Psi^r_s(t)+t-1)$. Thus
$$
\{(\la_s^i(-1/2+ \la_s^r\, (-1/2)),\Psi_s^i(\Psi^r_s(t)+t-1))\}
\times [-1,1]\subset W^u(P,f_{t,s}).
$$
If $\Psi_s^i(\Psi^r_s(t)+t-1)$ is in the interior
of $D_{t,s}^+$ we are done. Otherwise,
$\Psi_s^{i+1}(\Psi^r_s(t)+t-1)$ is the right extreme of
$D_{t,s}^+$. In this case the result follows using the $\la$-lemma:
$W^u(P,f_{t,s})$ accumulates 
the previous point at the right, 
recall the proof of Lemma~\ref{l.Wuptsegtovert}.
\end{demo}

\begin{lemm}
\label{l.semnome}
For every $x\in [-1,1]$, 
 $W^s(Q,f_{t,s})$
meets transversely
$\{x\}\times D_{t,s}^+\times [-1,1]$.  
\end{lemm}

\begin{demo}
Observe that, by construction, 
the point $(0,0,-1/2)$ belongs to $W^u(Q,f_{t,s})$.
Define $j>0$ by
$\Psi_s^{-j}(0)\in D_{t,s}^+$
(where $j\to \infty$ as $s\to 0$).
It is now immediate that 
$$
H=[-1,1]\times\{(\Psi_s^{-j}(0),-\la_u^{-j}\, (1/2))\}
\subset W^s(Q,f_{t,s})\quad
\mbox{and}
\quad
 H\pfrk (\{x\}\times D_{t,s}^+\times [-1,1])
\ne \emptyset,
$$
ending the proof of the lemma.
\end{demo}

We are now ready to prove  
Proposition~\ref{p.faixavert2}.
By the $\la$-lemma and Lemma~\ref{l.segto},
the forward orbit of the disk $\chi$ 
contains a strip $\De$ of the form
$\{x\}\times [a,b]\times [-1,1]$, $x\in [-1,1]$,
where $\alpha_0=[a,b]\subset D_{t,s}^+$.
Using the map $\cR_{t,s}$ and exactly arguing as in 
Section~\ref{ss.homoclinicclasses},
we inductively define strips 
$\De_k=\{x_k\}\times \alpha_k\times [-1,1]$, $x_k\in [-1,1]$
and $\alpha_k\subset D_{t,s}^+$,
such that
$\De_{k+1}\subset f_{t,s}^{n_k}(\De_k)$ and
$\alpha_{k+1}=\cR_{t,s}(\alpha_k)$.
By Lemma~\ref{l.J},
there is a first $k\in \NN$ such that $\alpha_{k+1}=\cR_{t,s}(\alpha_k)$ contains
$D_{t,s}^+$.
By Lemma~\ref{l.semnome}, $W^s(Q,f_{t,s})\pfrk \De_{k+1}\ne\emptyset$, 
thus   $W^s(Q,f_{t,s})\pfrk \De\ne\emptyset$,
ending the proof of the proposition.

\subsection{End of the proof of Theorem~\ref{t.main}}
\label{ss.endofthmain}
We now construct a one-parameter family of
diffeomorphisms 
$(g_s)$ satisfying Theorem~\ref{t.main}.
For that consider the arc $f_{t,s}$ defined as in the beggining of
 Section~\ref{s.colapso}.
We fix small $\bar t>0$
and consider the arc $g_s=f_{\bar t,-s}$. The results in the
previous section imply that
\begin{itemize}
\item
for every $s<0$,
$H(P,g_{s})$ and
$H(Q,g_{s})$ are non-hyperbolic and disjoint
(Theorem~\ref{t.iguais}),
\item
$s=0$,
$H(P,g_{0})=H^+(S, g_{0})$ and
$H(Q,f_{0})=H^-(S, g_{0})$,
(Theorems~\ref{t.latlat} and \ref{t.HssubsetHp}, see the beginning of  
Section~\ref{ss.beforecollapse}),
\item
for every $s>0$,
$H(P,g_{s})=H(Q,g_{s})$, (Theorem~\ref{t.HP=HQ}).
\end{itemize}

To finish the proof of Theorem~\ref{t.main} we need to see that
$\{S\}\in H(P,g_{s})\cap H(Q,g_{s})$
and to describe the maximal invariant set of $g_s$ in
the neighborhood $W$.

We assume that the neighborhood of the cycle $W$
is a {\em level of a fitration\/} 
of $f_{0,0}$ (thus, by continuity and compacity,
it is also a level of
a filtration for $f_{t,s}$ for every small $s$ and $t$).
This  means that there are compact sets $M_2$ and $M_1$,
$M_1$ contained in the interior of $M_2$, such that
$M_2\setminus M_1 = W$ and
$f_{0,0}(M_i)$ is contained in the interior of $M_i$, 
$i=1,2$. 
Hence, if $x\in W$ and 
$f^i_{t,s}(x)\not\in W$ for some $i$, then
$x$ is wandering: suppose, for instance,  
that $f^{i_0}_{t,s}(x)\in M_1$, where
$i_0>0 $. Then, there is a neighborhood $U_x\subset W$ of $x$ 
such that   $f^{i_0}_{t,s}(U_x)\subset M_1$. The definition of the
filtration implies that  $f^{i_0+j}_{t,s}(U_x)\subset M_1$ for all 
$j\ge 0$. 
Thus 
 $f^{i_0+j}_{t,s}(U_x)\cap U_x=\emptyset$ for all $j\ge 0$.
By shrinking $U_x$, we can assume that
$U_x, f_{t,s}(U_x), \dots, f_{t,s}^{i_0}(U_x)$ are pairwise 
disjoint. Therefore 
 $f^{j}_{t,s}(U_x)\cap U_x=\emptyset$ for all $j>0$, and $x$ is
wandering.

\medskip

Using the definition of the arc $f_{t,s}$,
it is immediate to check  the following:
let $\La_{t,s}$ be the maximal invariant set
of $f_{t,s}$ in $W$.

\begin{rema}
\label{r.s0}
For small positive $s$,
\begin{itemize}
\item
Every point
$(x,y,s)\in C=[-1,1]\times [-2,2]\times [-1,1]$
with $y\in (-\sqrt{s},\sqrt{s})$ 
(resp, $y\in [-2,-1)$ or $y\in (1,2)$)
is wandering.
\item
Consider $w=(x,y,z)\in C\cap \La_{t,s}$ and, 
for $i\in \ZZ$, let $w_i=g_s^i(w)=f_{\bar t,s}^i(w)$. If 
$w_i\in C$
we let 
$w_i=(x_i,y_i,z_i)$. 
Suppose that
$y_i\in [\sqrt{s},1]$ and $y_j\in [-2,\sqrt{s})$ for some
$j>0$. Then, for every $w_n\in C$
with $n\ge j$, one has $y_n\in [-2,\sqrt{s})$.
\end{itemize}
For $s=0$,
\begin{itemize}
\item
Every point
$(x,y,s)\in C$
with  $y\in [-2,-1)$ or $y\in (1,2)$
is wandering.
\item
Consider $w=(x,y,z)\in C\cap \La_{t,0}$ and, 
for $i\in \ZZ$, let $w_i=g_0^i(w)=f_{\bar t,0}^i(w)$. If 
$w_i\in C$
we let 
$w_i=(x_i,y_i,z_i)$. 
Suppose that
$y_i\in (0,1]$ and $y_j\in [-2,0)$ for some
$j>0$. Then, for every $w_n\in C$
with $n\ge j$, one has $y_n\in [-2,0)$.
\end{itemize}
\end{rema}
 
The  remark immediately
implies that, for small $s>0$,
$$
\La_{\bar t,s}\cap \Om(f_{\bar t,s})=\La_{\bar t,s}^+ \cup \La_{\bar t,s}^-,
$$
where
$\La_{\bar t,s}^+$ (resp. $\La_{\bar t,s}^-$)  
is the set of points 
$w\in \La_{\bar t,s}\cap \Om(f_{\bar t,s})$
such that, for every $i\in \ZZ$ with $w_i=(x_i,y_i,z_i)\in C$ 
one has $y_i\in [\sqrt{s},1]$
(resp. $y_i\in [-1,\sqrt{s}]$).
As in Section~\ref{s.pontosnaoerrantes},
$
H(P,f_{\bar t,s})=\La_{\bar t,s}^+$ and 
$H(Q,f_{\bar t,s})=\La_{\bar t,s}^-$.
Observe that we need to exclude the segment
$\{0\}\times (-\sqrt{s},\sqrt{s})\times \{0\}$ contained
in $\La_{\bar t,s}$ and consisting of wandering points.

For the saddle-node parameter $s=0$, one argues similarly. 
First, 
as above, one has that
$$
\La_{\bar t,0}=
(\La_{\bar t,0}\cap \Om(f_{\bar t,0}))=
\La_{\bar t,0}^+ \cup \La_{\bar t,0}^-,
$$
where
$\La_{\bar t,0}^+$
(resp. $\La_{\bar t,0}^-$)
 is the set of points 
$w\in \La_{\bar t,0}$,
such that, for every $i\in \ZZ$ with  
$w_i=(x_i,y_i,z_i)\in C$, one has $y_i\in [0,1]$
(resp. $y_i\in [-1,0]$).
As in the case $s>0$,
$H(P,f_{\bar t,0})=\La_{\bar t,0}^+$ and 
$H(Q,f_{\bar t,0})=\La_{\bar t,0}^-$.
It is immediate that $\La_{\bar t,0}^+\cap \La_{\bar t,0}^-=\{S\}$.

For parameters $s<0$ the result follows similarly (but now
the situation is much more simple).

\subsection{Proof of Theorem~\ref{t.saturated}}
Clearly, the homoclinic classes 
$H(P,f_{\bar t,0})$ and 
$H(Q,f_{\bar t,0})$ are not saturated.
We claim that there is not any transitive saturated
set $\Si$ containing $H(P,f_{\bar t,0})$. Otherwise,
the set $\Si$ must also contain $H(Q,f_{\bar t,0})$.
Thus $\Si$ contains the whole $\La_{\bar t,0}$. Using the filtration one has
$\Si=\La_{\bar t,0}$. But this set is not transitive:
Remark~\ref{r.s0} implies that there is no orbit 
going from a small neighborhood $U_Q$ of  $Q$
to   a small neighborhood $U_P$ of  $P$
and thereafter reuturning to $U_Q$.
This contradiction ends the proof of the theorem.

\medskip

\begin{tabular*}{220mm}{l l%
@{\extracolsep{4pt}}}
{\bf Lorenzo J. D\'\i az}  &  {\bf Bianca Santoro} \\

(lodiaz@mat.puc-rio.br) &  (bsantoro@math.mit.edu) \\

 Departamento de Matem\'{a}tica   & Department of Mathematics \\

PUC-Rio &   Massachusetts Institute of Technology \\

Marqu\^{e}s de S. Vicente 225  & 77 Massachusetts Avenue\\
22453-900 Rio de Janeiro RJ   &Cambridge MA 02139 \\
Brazil &  USA

\end{tabular*}


\medskip

\begin{thebibliography}{ZZZZZ}

\bibitem[Ab]{Ab} F. Abdenur,
{\em Generic robustness of spectral decompositions,\/}
to appear in Ann. Scient. Ecole Nor. Sup. Paris.



\bibitem[Ar]{Ar} M.C. Arnaud,
{\em Creation de conections en topologie $C^1$,\/}
Ergodic Th. \& Dyn. Systems, {\bf 21}, 339-381, (2001).

\bibitem[BD$_1$]{BD1} Ch. Bonatti e L.J. D\'\i az,
{\em Persistence of transitive diffeomorphims,\/}
Ann. Math., {\bf 143}, 367-396,  (1995).

\bibitem[BD$_2$]{BD01}
Ch. Bonatti and L.J. D\'\i az,
{\em On maximal transitive sets of generic diffeomorphisms,\/}
preprint PUC-Rio, 2001.


\bibitem[BDP]{BDP} Ch. Bonatti, L.J. D\'\i az, and E.R. Pujals, {\em A $\cC^1$-generic dichotomy for diffeomorphisms: Weak forms of hyperbolicicity or or infinitely many sinks or sources,\/}
to appear in Annals of Math.

\bibitem[CMP]{CMP}
C. M. Carballo, C. A. Morales e M.J. Pacifico,
{\em Homoclinic classes for generic $C^1$ vector fields,\/}
to appear in Ergodic Th. \& Dyn. Systems.



\bibitem[D$_1$]{Detds}
L.J. D\'\i az,
{\em
Robust nonhyperbolic dynamics e heterodimensional cycles,\/}
Ergodic Th. \& Dyn. Systems, {\bf 15},  291-315, (1995).

\bibitem[D$_2$]{Dnon}
L.J. D\'\i az,
{\em Persistence of  cycles and nonhyperbolic dynamics at the unfolding
of heteroclinic bifurcations,\/}
 Nonlinearity, {\bf 8}, 693-715, (1995).


\bibitem[DR$_1$]{DRnc}
L.J. D\'\i az and  J. Rocha,
{\em Non-conected heterodimensional cycles: bifurcation and stability\/}
Nonlinearity, {\bf 5},
1315-1341, (1992).


\bibitem[DR$_2$]{DRnon}
L.J. D\'\i az and  J. Rocha,
{\em Large measure of hyperbolic dynamics when unfolding
 heteroclinic cycles,\/}  Nonlinearity, {\bf 10},
857-884, (1997).


\bibitem[DR$_3$]{DRdss}
L.J. D\'\i az and J. Rocha,
{\em Noncritical saddle-node cycles and robust nonhyperbolic
dynamics,\/}  Journal of Dynamics and Stability
of Systems, {\bf 12}(2), 109-135, (1997).


\bibitem[DR$_4$]{DRetds}
L.J. D\'\i az and  J. Rocha,
{\em Partially hyperbolic and transitive  dynamics generated by
 heteroclinic cycles,\/}
 Ergodic Th. \& Dyn. Systems, {\bf 21}, 25-76,(2001).

\bibitem[DR$_5$]{DRfund}
L.J. D\'\i az and J. Rocha,
{\em Heterodimensional cycles, partial hyperbolicity
and limit dynamics,\/}
 to appear in Fundamenta Mathematicae.

\bibitem[DR$_6$]{DRprep}
L.J. D\'\i az and  J. Rocha,
{\em
Modelling collisions of hyperbolic homoclinic classes,\/}
in preparation.

\bibitem[DU]{DU}
L.J. D\'\i az and  R. Ures,
{\em Persistent homoclinic tangencies at the unfolding of cycles,\/}
Ann. Inst. Henri Poincar\'e, Analyse non lin\'eaire, {\bf 11}, 643-659,
(1994).


\bibitem[HPS]{HPS}
M. Hirsch, C. Pugh, et M. Shub,
{\em Invariant manifolds,\/}
Lecture Notes in Math.,
{\bf 583},
Springer Verlag,
(1977).

\bibitem[Ma]{Ma}
J. Mather,
{\em Commutators of diffeomorphisms,\/}
Comm. Math. Helv.,
{\bf 48},
195--233,
(1973).


\bibitem[Sm]{Sm} S. Smale,
{\em Differentiable dynamical systems,\/}
Bull. A.M.S.,
{\bf 73},
147-817,
(1967).

\end{thebibliography}
\end{document}